\documentclass{article}
\usepackage[letterpaper, portrait, margin=2.5cm]{geometry}

\usepackage{amsmath,amssymb,amsthm}
\usepackage{mathtools}
\usepackage{xcolor,colortbl}
\usepackage[all,cmtip]{xy}
\input diagxy

\usepackage{hyperref}

\usepackage{todonotes}

\newcommand{\catname}[1]{{\textsf{\textbf{#1}}}}

\newcommand{\Fieldsk}{\catname{Fields/k}}
\newcommand{\MSETS}{\catname{MarkedSets}}
\newcommand{\Sets}{\catname{Sets}}
\newcommand{\Cat}{\catname{Cat}}
\newcommand{\Varn}{\catname{Var}^\nu}
\newcommand{\op}{op}

\newcommand{\CC}{\mathbb C}

\newcommand{\FF}{\mathbb F}

\newcommand{\PP}{\mathbb P}

\newcommand{\VV}{\mathbb V}

\newcommand{\ZZ}{\mathbb Z}

\newcommand{\mce}{\mathcal E}

\newcommand{\lbm}{\left[ \begin{matrix}}
\newcommand{\rem}{\end{matrix} \right]}
\newcommand{\lp}{\left(}
\newcommand{\rp}{\right)}
\newcommand{\lb}{\left\{}
\newcommand{\rb}{\right\}}
\newcommand{\lav}{\left|}
\newcommand{\rav}{\right|}

\newcommand{\SL}{\sum\limits}

\newcommand{\ra}{\rightarrow}
\newcommand{\bra}{\dashrightarrow}

\newcommand{\RA}{\Rightarrow}

\newcommand{\racts}{\curvearrowright}

\newcommand{\injects}{\hookrightarrow}

\DeclareMathOperator{\PGL}{PGL}

\DeclareMathOperator{\PSL}{PSL}

\newcommand{\st}{ \ | \ }

\newcommand{\wt}{\widetilde}

\DeclareMathOperator{\Aut}{Aut}

\DeclareMathOperator{\im}{Im}

\DeclareMathOperator{\spec}{Spec}

\DeclareMathOperator{\Bran}{Bran}
\DeclareMathOperator{\Fin}{Fin}

\DeclareMathOperator{\Gal}{Gal}
\DeclareMathOperator{\trdeg}{tr.deg}

\DeclareMathOperator{\Sol}{Sol}
\DeclareMathOperator{\Ab}{Ab}

\DeclareMathOperator{\Sym}{Sym}

\DeclareMathOperator{\RD}{RD}
\DeclareMathOperator{\ed}{ed}
\DeclareMathOperator{\RDC}{\operatorname{RD}_{\mathbb{C}}}
\newcommand{\KD}{K^{(d)}}
\newcommand{\KE}{K^{\mathcal{E}}}
\newcommand{\KS}{K^{\text{Sol}}}
\newcommand{\KSEP}{K^{\text{sep}}}

\DeclareMathOperator{\MOO}{M_{11}}
\DeclareMathOperator{\MOW}{M_{12}}
\DeclareMathOperator{\MWW}{M_{22}}
\DeclareMathOperator{\MWH}{M_{23}}
\DeclareMathOperator{\MWF}{M_{24}}
\DeclareMathOperator{\COO}{Co_1}
\DeclareMathOperator{\COW}{Co_2}
\DeclareMathOperator{\COH}{Co_3}
\DeclareMathOperator{\JO}{J_1}
\DeclareMathOperator{\JW}{J_2}
\DeclareMathOperator{\JH}{J_3}
\DeclareMathOperator{\JF}{J_4}
\DeclareMathOperator{\SUZ}{Suz}
\DeclareMathOperator{\MCL}{McL}
\DeclareMathOperator{\HS}{HS}
\DeclareMathOperator{\FIWW}{Fi_{22}}
\DeclareMathOperator{\FIWH}{Fi_{23}}
\DeclareMathOperator{\FIWF}{Fi_{24}'}
\DeclareMathOperator{\THO}{Th}
\DeclareMathOperator{\HN}{HN}
\DeclareMathOperator{\HE}{He}
\DeclareMathOperator{\B}{B}
\DeclareMathOperator{\MO}{M}
\DeclareMathOperator{\ON}{O'N}
\DeclareMathOperator{\RU}{Ru}
\DeclareMathOperator{\LY}{Ly}
\newcommand{\PERM}{\text{Perm}}

\DeclareMathOperator{\SCH}{Sch}
\newcommand{\GAP}{\texttt{GAP}}
\newcommand{\arxiv}{\texttt{arXiv}}
\newcommand{\SAGE}{\texttt{SageMath}}
\newcommand{\ATLAS}{\mathbb{A}\mathbb{T}\mathbb{L}\mathbb{A}\mathbb{S}}

\newcommand{\TX}{\prescript{T}{}{X}}

\newcommand{\TXG}{\prescript{T}{}{X_G}}
\newcommand{\TYG}{\prescript{T}{}{Y_G}}
\newcommand{\TZG}{\prescript{T}{}{Z_G}}
\DeclareMathOperator{\SPAN}{Span}

\theoremstyle{definition}
\newtheorem{definition}{Definition}[section]
\newtheorem{problem}[definition]{Problem}

\newtheorem{example}[definition]{Example}
\newtheorem{notation}[definition]{Notation}

\theoremstyle{remark}
\newtheorem{remark}[definition]{Remark}

\theoremstyle{plain}
\newtheorem{prop}[definition]{Proposition}
\newtheorem{lemma}[definition]{Lemma}
\newtheorem{theorem}[definition]{Theorem}
\newtheorem{corollary}[definition]{Corollary}
\newtheorem{conjecture}[definition]{Conjecture}

\title{Generalized Versality, Special Points, and Resolvent Degree for the Sporadic Groups}
\author{Claudio G\'{o}mez-Gonz\'{a}les, Alexander J. Sutherland, and Jesse Wolfson\footnote{The third author was supported in part by NSF Grant DMS-1944862.}}
\date{}

\begin{document}

\maketitle

\begin{abstract}
Resolvent degree is an invariant measuring the complexity of algebraic and geometric phenomena, including the complexity of finite groups. To date, the resolvent degree of a finite simple group $G$ has only been investigated when $G$ is a cylic group; an alternating group; a simple factor of a Weyl group of type $E_6$, $E_7$, or $E_8$; or $\PSL\lp 2, \FF_7 \rp$. In this paper, we establish upper bounds on the resolvent degrees of the sporadic groups by using the invariant theory of their projective representations. To do so, we introduce the notion of (weak) $\RD_k^{\leq d}$-versality, which we connect to the existence of ``special points'' on varieties.
\end{abstract}

\setcounter{tocdepth}{3}
\tableofcontents

\section{Introduction}\label{sec:Introduction}

What is the least $d$ for which a solution of the general degree $n$ polynomial admits a formula using only (algebraic) functions of $d$ or fewer variables?  As Abel realized, the general degree $n$ polynomial has Galois group $S_n$ and this question, in modern language, asks for the \textbf{resolvent degree} of the symmetric group, denoted
\begin{equation*}
\RDC(S_n)=\RDC(A_n)=:\RDC(n),
\end{equation*}
an invariant first introduced independently by Brauer \cite{Brauer1975} and Arnol'd-Shimura \cite{ArnoldShimura1976}.  To the best of our knowledge, Klein was the first to consider this question for other finite groups, most notably the group $\PSL\lp 2, \FF_7 \rp$ \cite{Klein1879}. Note that for a finite group $G$ with Jordan-H\"{o}lder decomposition $\lb G_1,\dotsc,G_s \rb$, \cite[Theorem 3.3]{FarbWolfson2019} yields that
\begin{equation*}
	\RD_k(G) \leq \max\lb \RD_k(G_1), \dotsc, \RD_k(G_s) \rb,
\end{equation*}
\noindent
with equality if every $G_j$ can be realized as a subgroup of $G$.\footnote{In fact, it seems reasonable to expect that the inequality can be made into an equality without assumptions on $G$.}

Following Klein and the classification of finite simple groups, one is led to the following question \cite[Problem 3.5]{FarbWolfson2019}:
\begin{problem}[$\RD_k(G)$ for Finite Simple $G$]\label{prob: RD of Finite Simple Groups}
Compute the resolvent degree of all finite simple groups $G$.
\end{problem}

To date, Problem \ref{prob: RD of Finite Simple Groups} has only been addressed by providing upper bounds on $\RD_k(G)$ when $G$ is a cyclic group (in which case $\RD_k(G) \equiv 1$),
an alternating group (see \cite{HeberleSutherland2023,Sutherland2021,Wolfson2021}, or \cite{Hamilton1836,Hilbert1927,Klein1879,Klein1884,Klein1887,Klein1905,Segre1945,Tschirnhaus1683} for classical references), when $G$ is a simple factor of a Weyl group of type $E_6$, $E_7$, or $E_8$ \cite{FarbWolfson2019,FarbKisinWolfson2023,Reichstein2022}, or when $G=\PSL\lp 2,\FF_7 \rp$ (see \cite{Klein1879} for the classical reference or \cite{FarbKisinWolfson2023} for a modern version). 

The Classification of Finite Simple Groups says that a finite simple group $G$ falls into one of four categories:
\begin{enumerate}
\item $G$ is cyclic of prime order;
\item $G$ is an alternating group $\lp A_n, n \geq 5 \rp$;
\item $G$ is a simple group of Lie type (of which there are 16 families); or
\item $G$ is one of 26 finite simple groups that do not belong to one of the infinite families above.
\end{enumerate}

\noindent
The 26 groups in (4) are known as the \textbf{sporadic groups}. In this paper, we investigate Problem \ref{prob: RD of Finite Simple Groups} by giving upper bounds on $\RD_k(G)$ for all sporadic groups $G$. For each group $G$, we use the invariant theory of a projective representation over $\CC$ of minimal dimension to construct a complex $G$-variety $X_G$ with the property that $\RD_k(G)\le \dim_{\CC}(X_G)$.  For $G$ one of the Mathieu groups $\MOO,\MOW,\MWH,\MWF$, we prove nothing new: the projective representation in these cases is just the projectivization of the permutation representation, $X_G$ is the vanishing locus of the first four ($\MOO,\MOW$) or five ($\MWH,\MWF$) elementary symmetric polynomials, and the bounds on $\RD_k(G)$ follow from the bounds on $\RD_k(S_n)$ which appear previously in the literature (with $S_{11}$ and $S_{12}$ in \cite{Segre1945}, and $S_{23}$ and $S_{24}$ in \cite{Sutherland2021}). Our primary interest is thus the remaining 22 sporadic groups.  Here, our bounds appear to be genuinely new, and we obtain them by proving that the variety $X_G$ is ``$\RD_{\CC}^{\le d}$-versal'' for $d<\dim_{\CC}(X_G)$ (see Definition \ref{def:RDKD versality}, Lemma \ref{lem:RDKD-Versality via $G$-Equivariant Rational Correspondences}, and Theorem \ref{thm:Bounds on the Resolvent Degree of the Sporadic Groups} for the construction of $X_G$), from which the bound on $\RD_k(G)$ follows by \cite[Proposition 3.10]{FarbWolfson2019}. More explicitly, we have:

\begin{corollary}[Appears as Corollary \ref{cor:Explicit Form of Bounds on the Resolvent Degree of the Sporadic Groups}: Explicit Form of Theorem \ref{thm:Bounds on the Resolvent Degree of the Sporadic Groups}]
For any field $k$, we have
\begin{align*}
	\RD_k(\JW) &\leq 5, &\RD_k(\MWF) &\leq 18, &\RD_k(\HE) &\leq 48, &\RD_k(\FIWH) &\leq 776,\\
	\RD_k(\MOO) &\leq 6, &\RD_k(\HS) &\leq 18, &\RD_k(\JO) &\leq 51, &\RD_k(\FIWF) &\leq 779,\\
	\RD_k(\MOW) &\leq 7, &\RD_k(\MCL) &\leq 19, &\RD_k(\FIWW) &\leq 74, &\RD_k(\JF) &\leq 1328,\\
	\RD_k(\MWW) &\leq 8, &\RD_k(\COH) &\leq 20, &\RD_k(\HN) &\leq 129, &\RD_k(\LY) &\leq 2475,\\
	\RD_k(\SUZ) &\leq 10, &\RD_k(\COW) &\leq 20, &\RD_k(\THO) &\leq 244, &\RD_k(\B) &\leq 4365,\\
	\RD_k(\JH) &\leq 16, &\RD_k(\COO) &\leq 21, &\RD_k(\ON) &\leq 338, &\RD_k(\MO) &\leq 196874.\\
	\RD_k(\MWH) &\leq 17, &\RD_k(\RU) &\leq 26, 
\end{align*}
\end{corollary}

The proof of Theorem \ref{thm:Bounds on the Resolvent Degree of the Sporadic Groups} is comprised of two distinct parts, which are set up in Sections \ref{sec:Generalized Versality and Special Points} and \ref{sec:Resolvent Degree and the Sporadic Groups}, respectively. First, in \cite[Definition 3.8]{FarbWolfson2019}, notions of solvable versality and $\RD_k$-versality for finite groups $G$ were introduced. In \cite{FarbKisinWolfson2023}, this was generalized by defining 1) a notion of a {\em class of accessory irrationalities} $\mathcal{E}$ \cite[Definition 4.1]{FarbKisinWolfson2023} and, given such a class $\mathcal{E}$, by defining 2) a notion of {\em $\mathcal{E}$-versality} \cite[Definition 4.4]{FarbKisinWolfson2023}.  A connection to resolvent degree was briefly discussed in \cite[Paragraph 4.1.3 and Lemma 4.9]{FarbKisinWolfson2023}, along with a discussion of historical roots of these notions of generalized versality and a call to better understand them \cite[Remark 4.10]{FarbKisinWolfson2023}.  We build on this framework here, introducing new examples of classes of accessory irrationalities $\mathcal{E}=\RD_k^{\le d}$ (``accessory irrationalities of resolvent degree at most $d$'') which have been implicit in the literature, and broadening \cite[Definition 4.4]{FarbKisinWolfson2023} and the attendant lemmas to allow for arbitrary algebraic groups $G$. 

Second, Duncan and Reichstein \cite[Theorem 1.1 (a,b)]{DuncanReichstein2015} reframed versality in terms of rational points on twisted forms. Building on this, we connect the above generalizations of versality to the existence of ``special points''---a catch-all term we use to encompass points defined via any specified class $\mathcal{E}$ of accessory irrationalities---as in Theorem \ref{thm:Introduction Theorem}.  This perspective, with $\mathcal{E}=\RD_k^{\le d}$, is implicit in \cite[Lemma 14.5]{Reichstein2022}, and our hope is that by making it explicit, we can make it more widely known and used.  Concretely, for a field $K$ and a class of accessory irrationalities $\mathcal{E}$, let $K^{\mathcal{E}}$ denote an $\mathcal{E}$-closure of $K$ (see Section \ref{subsec:Accessory Irrationalities}).

\begin{theorem}[Appears as Theorem \ref{thm:RDKD-versality and Special Points}: Generalized Versality and Special Points]\label{thm:Introduction Theorem}
Let $G$ be an algebraic group over $k$ and let $X$ be an irreducible, generically free, quasiprojective $G$-variety. Let $\mathcal{E}$ be a class of accessory irrationalities.  Then: 
\begin{enumerate}
\item $X$ is weakly $\mathcal{E}$-versal if and only if for every $G$-torsor $T\ra \spec(K)$ with $K$ finitely generated over $k$, $\TX\lp K^{\mathcal{E}}\rp\not= \emptyset$.
\item If $G$ is smooth, $X$ is $\mathcal{E}$-versal if and only if for every $G$-torsor $T \ra \spec(K)$ with $K$ finitely generated over $k$, $K^{\mathcal{E}}$-points are dense in $\TX$. 
\end{enumerate}
\end{theorem}

In Section \ref{sec:Resolvent Degree and the Sporadic Groups}, we apply the framework of $\mathcal{E}$-versality outlined in Section \ref{sec:Generalized Versality and Special Points} to the 22 sporadic groups not equal to $\MOO,\MOW,\MWH,$ or $\MWF$. More specifically, we use the invariant theory of a projective representation for each sporadic group $G$ to construct a complex variety $X_G$ which we then show is $\RD_{\CC}^{\leq d_G}$-versal for some $d_G \le \dim_{\CC}\lp X_G\rp$. Results of Reichstein, namely \cite[Theorems 1.2, 1.3]{Reichstein2022}, then allow us to conclude the upper bounds on $\RD_k(G)$ for all fields $k$.

\paragraph{Outline of the Paper}
The remainder of this paper proceeds as follows. In Section \ref{sec:Background}, we introduce the relevant background on torsors (Section \ref{subsec:Torsors}), accessory irrationalities (Section~\ref{subsec:Accessory Irrationalities}) and resolvent degree (Section \ref{subsec:Resolvent Degree}). In Section \ref{sec:Generalized Versality and Special Points}, we recall the framework of $\mathcal{E}$-versality and connect it to the existence of $\mathcal{E}$-points. In Section \ref{sec:Resolvent Degree and the Sporadic Groups}, we establish upper bounds on the resolvent degree of the sporadic groups (Section \ref{subsec:Upper Bounds on the Resolvent Degree of the Sporadic Groups}) and understand these bounds in terms of the prior literature (Section \ref{subsec:Context for Sporadic Group Bounds}).

\paragraph{Supplementary Materials}
Some results in this paper rely on calculations performed with \GAP \ \cite{GAP2022} and \SAGE \ \cite{SageMath2022}. Along with this paper, we have included supplementary files which can be used to verify all computations. They can be found on the \arxiv \ submission of this paper, as well as at the following web address:

\begin{center}
	\href{https://www.alexandersutherland.com/research/RD-for-the-Sporadic-Groups}{\texttt{https://www.alexandersutherland.com/research/RD-for-the-Sporadic-Groups}} \ .
\end{center}

\noindent
In addition to the \SAGE \ script, we have also included the relevant data and outputs as plain text files. The supplementary files on the \arxiv \ can be accessed by downloading the source package.

\paragraph{Correctness of Information for Sporadic Groups}
The computations for the Molien series (or the coefficients thereof) for projective representations of sporadic groups in this project rely on the \GAP \ character table library \cite{Breuer2013} of the computer algebra system \GAP \ \cite{GAP2022}. We refer the reader to \cite{BreuerMalleOBrien2017} for details on correctness, while noting the excerpt ``all character tables contained in the $\ATLAS$, incorporating the corrections, and many more, are stored electronically in the character table library \cite{Breuer2013} of the computer algebra system \GAP \ \cite{GAP2015}.''

\paragraph{Conventions}
\begin{enumerate}
\item We define a $K$-variety to be a quasiprojective scheme of finite type over $K$. We do not require varieties to be reduced or irreducible.
\item For a collection of homogeneous polynomials $\lb f_1,\dotsc,f_s \rb \subseteq K\left[ x_0,\dotsc,x_n \right]$, we write $\VV\lp f_1,\dotsc,f_s \rp$ for the subvariety of $\PP_K^n$ determined by the (scheme-theoretic) intersection $f_1 = \cdots = f_s = 0$.
\item We denote the $K$-points of a variety $X$ by $X(K)$.
\item Unless otherwise specified, when we refer to an algebraic group, we mean an arbitrary algebraic group (it need not be finite, linear, nor smooth).
\item Given a field $K$, we denote a separable closure of $K$ by $\KSEP$. 
\item For maps between varieties, we denote regular morphisms by $\ra$, and rational maps by $\bra$.
\end{enumerate}

\paragraph{Acknowledgements}
The authors would like to thank Aaron Landesman and Daniel Litt for helpful comments on a draft, and Alice Silverberg for help with a reference. The authors thank the anonymous referee for helpful comments and suggestions.


\section{Background}\label{sec:Background}

In this section, we recall the necessary background on torsors, accessory irrationalities, and resolvent degree. We fix a ground field $k$.   

\subsection{Torsors}\label{subsec:Torsors}
Consider an algebraic group $G$ over $k$, i.e. a group scheme of finite type over $k$.  We refer the reader to \cite[Chapter 5]{Poonen2017} for a good summary and list of references for algebraic groups.  For example, recall that results of Chow and Conrad (see e.g. \cite[Theorem 5.2.20]{Poonen2017}) show that $G$ is automatically a quasi-projective $k$-scheme.  Let $X$ be a quasiprojective $G$-variety over $k$. In this section, $K$ always denotes a finitely generated $k$-field.

\begin{definition}[$G$-Torsors]\label{def:G-Torsors}
A \textbf{right} (respectively, left) \textbf{$G$-torsor over $X$} is a flat morphism $Y \ra X$ of $k$-schemes such that $G$ acts on $Y$ on the right (respectively, left) by $\sigma:G \times Y \ra Y$ and such that the map
\begin{align*}
	G \times Y &\ra Y \times_X Y\\
	(g,y) &\mapsto (\sigma(g,y),y),
\end{align*}

\noindent
is an isomorphism. We say that the $G$-torsor $Y \ra X$ is \textbf{split} if it admits a section; this is equivalent to its class $[Y]$ in $H^1(k,G)$ being trivial.
\end{definition}

Next, we introduce twisted varieties.

\begin{definition}[Twists]\label{def:Twists}
Let $G$ be an algebraic group over $k$. Let $T \ra \spec(K)$ be a $G$-torsor. $G$ acts on $T \times X$ diagonally and yields a $G$-torsor $T \times X \ra \TX$. We say $\TX$ is \textbf{the twist of $X$ by $T$}.
\end{definition}

\noindent
Note that our assumption of quasiprojectivity of $X$ implies that $\TX$ is well-defined as the geometric quotient of $T \times X$ by $G$; see \cite[Section 5.12.5]{Poonen2017} for details when $G$ is smooth, and \cite[Proposition 2.12]{Florence2008} in general. 

In Section \ref{sec:Generalized Versality and Special Points}, we will want to move interchangeably between torsors over finitely generated $k$-fields and their integral models.

\begin{definition}[Integral Models]\label{def:Integral Model}
Given a $G$-torsor $T \ra \spec(K)$ over $k$, an \textbf{integral model} of $T \ra \spec(K)$ is a morphism $Y \ra Y/G$, where $Y$ is a generically free $G$-variety with $k\lp Y/G \rp = K$, along with an isomorphism $Y \times_{Y/G} \spec(K) \cong T$.\footnote{Note that our assumption that $k(Y/G)=K$ implies that $Y/G$ is irreducible. In particular, all integral models are birational to each other.}
\end{definition}

For further background on $G$-torsors, we refer the reader to \cite[Ch. 3, Sec. 4]{Milne1980} and \cite[Section 5.12]{Poonen2017} for general results and to \cite[Section 10]{Reichstein2022} for connections of $G$-torsors to essential dimension and resolvent degree.

\subsection{Accessory Irrationalities}\label{subsec:Accessory Irrationalities}
Accessory irrationalities appear prominently in work of Klein \cite{Klein1884} (see also \cite{Klein1879,Klein1887,Klein1905} and Chebotarev \cite{Chebotarev1932,Chebotarev1934}).  The first formal definition of which we are aware appears in \cite[Definition 4.1]{FarbKisinWolfson2023} in the language of branched covers.  We recall this here, introduce an equivalent formulation in the language of field extensions, give new examples of $\mathcal{E}$, and introduce the notion of a closure of a field with respect to a class of accessory irrationalities (Definition~\ref{def:RDKD closure}).

Following \cite[4.1.1]{FarbKisinWolfson2023}, we consider branched covers $p\colon Y\to X$ of normal $k$-varieties, i.e. dominant, finite maps. More generally, we will use {\bf branched cover} to refer to a generically finite, dominant rational map $p\colon Y\bra X$. Branched covers form a category in the usual fashion, and they are preserved under pullback in the following sense: if $f\colon X'\to X$ is any map of normal $k$-varieties, then we denote by $f^*Y$ the normalization of $Y\times_X X'$ and observe that $f^*p\colon f^*Y\to X'$ is again a branched cover.

\begin{definition}[Definition 4.1 of \cite{FarbKisinWolfson2023}]\label{d:aibra}
    Let $k$ be a field. Let $\Varn_k$ denote the category of normal $k$-varieties. Let 
    \[
        \Bran\colon(\Varn_k)^{\op}\to\Cat
    \]
    be the functor\footnote{We follow standard usage and do not distinguish between a functor and a pseudo-functor when the latter is the manifestly appropriate notion given the target in question.  Note that the natural isomorphisms required for pseudo-functoriality are canonical in this case, as they come from the universal property of normalization and fiber product.} which sends a finite type normal $k$-scheme $X$ to its category of branched covers. A {\bf class of accessory irrationalities} $\mathcal{E}$ is a subfunctor $\mathcal{E}\subset\Bran$ such that
    \begin{enumerate}
        \item For any $X$, $\mathcal{E}(X)\subset \Bran(X)$ is a full subcategory. 
        \item For any $X$, the identity $X\to X$ is in $\mathcal{E}(X)$.
        \item\label{it:coprod} $\mathcal{E}(X\coprod X')=\mathcal{E}(X)\times\mathcal{E}(X')$.
        \item If $E,E'\in\mathcal{E}(X)$, then $E\times_XE'\in\mathcal{E}(X)$.
        \item\label{it:accirrbirat} If $U\subset X$ is a dense open, then $\mathcal{E}(X)\to\mathcal{E}(U)$ is an equivalence of categories.
        \item If $E\to X'\to X$ are branched covers and if $E\to X$ is in $\mathcal{E}(X)$, then $E\to X'$ is in $\mathcal{E}(X')$.
    \end{enumerate}
\end{definition}

For the present paper, and to make explicit the connection to the perspective in \cite{Reichstein2022}, we rephrase this in the language of field extensions.

\begin{definition}[Definition 4.1 of \cite{FarbKisinWolfson2023} via Field Extensions]\label{def:Def 4.1 of FKW2022 via Field Extensions}
    Let $k$ be a field, let $\Fieldsk$ be the category of fields over $k$, and let
    \[
        \Fin\colon \Fieldsk\to\Cat
    \]
    be the functor which sends a $k$-field $K$ to the category of finite, semi-simple commutative $K$-algebras.\footnote{N.b. Wedderburn's theorem implies that any finite semi-simple commutative algebra over a field is a finite product of finite extensions of that field.} A {\bf class of accessory irrationalities} $\mathcal{E}$ is a subfunctor $\mathcal{E}\subset\Fin$ such that:
    \begin{enumerate}
        \item For all $K$, $\mathcal{E}(K)\subset \Fin(K)$ is a full subcategory.
        \item For all $K$, we have $K\in\mathcal{E}(K)$.
        \item If $E,E'\in \mathcal{E}(K)$, then $E\otimes_K E'\in\mathcal{E}(K)$.
        \item If $K\hookrightarrow L$ is a finite extension of $k$-fields, $L\hookrightarrow E$ is finite, and $K\hookrightarrow L\hookrightarrow E$ is in $\mathcal{E}(K)$, then $E\in \mathcal{E}(L)$.
    \end{enumerate}
\end{definition}

\begin{lemma}[Equivalence of Definitions]\label{l:aifin=aibra}
   Definitions \ref{d:aibra} and \ref{def:Def 4.1 of FKW2022 via Field Extensions} are equivalent: the assignment $X\mapsto k(X)$ induces an equivalence between the category of subfunctors of $\Bran$ satisfying the axioms of Definition~\ref{d:aibra} and the category of subfunctors of $\Fin$ satisfying the axioms of Definition~\ref{def:Def 4.1 of FKW2022 via Field Extensions}.
\end{lemma}

\begin{proof}
    Let $\mathcal{E}\subset\Bran$ be a class of accessory irrationalities. By assumptions~\ref{it:coprod} and~\ref{it:accirrbirat} of Definition~\ref{d:aibra}, we see that $\mathcal{E}$ is determined up to equivalence by its restriction to the sub-category of irreducible, affine, normal $k$-varieties. For any such $X$, we have a natural equivalence of categories
    \[
        \Fin\left(k(X)\right)\cong\{k(E)~|~E\to X \in\Bran(X)\}.
    \]
    Denote by $(k\circ \mathcal{E})(X)$ the full sub-category of $\Fin(k(X))$ determined by $\mathcal{E}(X)$ under the above equivalence. The assumptions of Definition~\ref{d:aibra} on $\mathcal{E}$ imply that $k\circ\mathcal{E}$ satisfies the assumptions of Definition~\ref{def:Def 4.1 of FKW2022 via Field Extensions}, as claimed. 

    It remains to show that any subfunctor $\mathcal{E}\subset\Fin$ satisfying Definition~\ref{def:Def 4.1 of FKW2022 via Field Extensions} arises as above.  Fix such a subfunctor $\mathcal{E}\subset \Fin$.  For $X\in \Varn_k$, let $\widetilde{\mathcal{E}}(X)\subset\Bran(X)$ be the full subcategory consisting of all branched covers $Y\dashrightarrow X$ such that for any irreducible component $X_i\subset X$, the restriction $Y|_{X_i}\dashrightarrow X_i$ has
    \[
        k(X_i)\hookrightarrow k(Y|_{X_i}) \in \mathcal{E}(k(X_i)).
    \]
    This definition ensures that $\widetilde{\mathcal{E}}\subset \Bran$ satisfies Assumptions~\ref{it:coprod} and~\ref{it:accirrbirat} of Definition~\ref{d:aibra}.  The remaining assumptions follow from the corresponding assumptions of Definition~\ref{def:Def 4.1 of FKW2022 via Field Extensions}, and by direct inspection, we see that $k\circ\widetilde{\mathcal{E}}\cong \mathcal{E}$ as claimed.
\end{proof}

Motivated by classical examples in the literature (see Example \ref{ex:Example 4.3(3) of FKW2022}), we introduce terminology for special classes of accessory irrationalities.

\begin{definition}[Saturation and Closure Under Extensions]\label{def:Saturation and Closure Under Extensions}
    Let $\mathcal{E}\colon \Fieldsk\to \Cat$ be a class of accessory irrationalities. 
    \begin{enumerate}
        \item We say $\mathcal{E}$ is {\bf saturated} if for all finite extensions of $k$-fields $K\hookrightarrow L$, that $K\hookrightarrow L\hookrightarrow E$ is in $\mathcal{E}(K)$ implies that $L\in \mathcal{E}(K)$.
        \item We say $\mathcal{E}$ is {\bf closed under extensions} if for all finite extensions of $k$-fields $K\hookrightarrow L$, $L\in\mathcal{E}(K)$ and $E\in\mathcal{E}(L)$ together imply that $K\hookrightarrow L\hookrightarrow E$ is in $\mathcal{E}(K)$.
    \end{enumerate}
\end{definition}

\begin{example}[Example 4.3(3) of \cite{FarbKisinWolfson2023}]\label{ex:Example 4.3(3) of FKW2022}\mbox{}
    \begin{enumerate}
        \item Let $\Ab(K)$ be the category of finite semisimple commutative $K$-algebras which split as products of abelian extensions of $K$. Then the assignment $K\mapsto\Ab(K)$ defines a saturated class of accessory irrationalities which is not closed under extensions. 
        \item Let $\Sol(K)$ be the full subcategory of finite semisimple commutative $K$-algebras which split as products of solvable extensions of $K$.  Then the assignment $K\mapsto \Sol(K)$ defines a saturated class of accessory irrationalities which is closed under extensions.
    \end{enumerate}

Note that $\Sol$ is the closure of $\Ab$ under extensions, i.e. it is the minimal class of accessory irrationalities which contains $\Ab$ and is closed under extensions.
\end{example}

\subsection{Resolvent Degree}\label{subsec:Resolvent Degree}
Resolvent degree was first defined independently by Brauer \cite{Brauer1975} and Arnol'd-Shimura \cite{ArnoldShimura1976} in the context of field extensions.  The first contemporary reference on resolvent degree is \cite{FarbWolfson2019}. We begin by reviewing the definition of resolvent degree for algebraic groups, and then introduce the notion of the resolvent degree of a functor.

\subsubsection{Resolvent Degree for Varieties, Fields, and Algebraic Groups}
We recall the definitions here and refer the reader to \cite{FarbWolfson2019,Wolfson2021,Reichstein2022} for more background. 

\begin{definition}[Essential Dimension of a Branched Cover of Varieties]\label{d:edbra}
    Let $Y\bra X$ be a branched cover of $k$-varieties (i.e. a generically finite, dominant rational map).  The \textbf{essential dimension} of $Y\bra X$ over $k$, denoted $\ed_k(Y\bra X)$, is the minimal $d$ for which there exists 
    \begin{enumerate}
        \item a branched cover $\tilde{Z}\to Z$ with $\dim_k Z=d$,
        \item a dense Zariski open $U\subset X$,
        \item a map $f\colon U\to Z$, and 
        \item an isomorphism $f^*\tilde{Z}\simeq Y|_U$ over $U$.
    \end{enumerate}
\end{definition}

\begin{definition}[Resolvent Degree of a Branched Cover of Varieties]\label{d:rdbra}
    Let $Y \bra X$ be a branched cover of $k$-varieties. The \textbf{resolvent degree} of $Y \bra X$ over $k$, denoted $\RD_k(Y \bra X)$, is the minimal $d$ for which there exists a tower of branched covers
    \begin{equation*}
        X_r \bra \cdots \bra X_0 = X,
    \end{equation*}

    \noindent
    with a factorization $X_r \bra Y \bra X$ such that $\ed_k(X_j \bra X_{j-1}) \leq d$ for all $1 \leq j \leq r$.
\end{definition}

The definition of the resolvent degree of an extension of $k$-fields extends straightforwardly to the case of finite, semisimple commutative algebras over $k$-fields, and the basic properties carry over as well. In particular, we have the following example.

\begin{example}[The class $\RD_k^{\le d}$]
    Fix $d\ge 0$. For a $k$-field $K$, let $\RD_k^{\le d}(K)$ denote the category of all finite semisimple commutative $K$-algebras $A$ such that $\RD_k(A/K)\le d$. The proof of \cite[Lemma 2.5(2)]{FarbWolfson2019} shows that the assignment $K\mapsto \RD_k^{\le d}(K)$ is indeed functorial, while that of \cite[Lemma 2.7]{FarbWolfson2019} shows that it satisfies the definition of a saturated class of accessory irrationalities which is closed under extensions.
\end{example}

Reichstein extended the above notion of essential dimension in \cite{Reichstein2000} as follows:
\begin{definition}[Essential Dimension of a $G$-Variety]
    Let $G$ be an algebraic group over $k$. Let $X$ be a generically free $G$-variety. The {\bf essential dimension} of $X\bra X/G$ is the least $d$ such that there exists a dominant, $G$-equivariant rational map $X\bra Y$ with $d=\dim(Y/G)$.
\end{definition}

Following \cite{Reichstein2021,Reichstein2022}, we build on this here.
\begin{definition}[Resolvent Degree of a $G$-Variety]\label{d:rdgvar}
    Let $G$ be an algebraic group. Let $X$ be a quasi-projective $G$-variety over $k$. The \textbf{resolvent degree} of $X\bra X/G$ is
    \begin{equation*}
        \RD_k(X\bra X/G) = \min\lb \max\{ \RD_k(E\bra X/G),\ed_k(X|_E\bra E)\} \rb.
    \end{equation*}
    where the minimum is over generically finite dominant maps $E\bra X/G$.
\end{definition}

\begin{remark}
    The above definition differs slightly from that of \cite{Reichstein2021,Reichstein2022}. We show in Lemma~\ref{l:rdgcomp} below that it agrees with Reichstein's.
\end{remark}

Recall that a $G$-variety $X$ is {\bf primitive} if $G$ acts transitively on the set of geometrically irreducible components of $X$.  A $G$-variety is {\bf generically free} if the locus of points with trivial (scheme theoretic) stabilizer is dense and open. We record the following elementary lemma for later use. 

\begin{lemma}[Faithful and Irreducible Implies Generically Free]\label{lem:faithful and irred suffices}
    Let $G$ be a finite group and $X$ an irreducible $G$-variety.  Then the action of $G$ on $X$ is generically free if and only if it is faithful.
\end{lemma}

\begin{proof}
    That generically free implies faithful is immediate. For the converse, note that given $g \in G \setminus \{1\}$, the (scheme theoretic) fixed set $X^g\subset X$ is a closed subset, which is not all of $X$ because the action is faithful. If $X$ is irreducible, then $X$ cannot be written as the union of a finite number of proper closed subsets. Therefore, $X\setminus \bigcup_{g\in G \setminus \{1\}} X^g$ is a non-empty open in which every point has trivial scheme theoretic stabilizer. As $X$ is irreducible, this open is dense, and thus the action is generically free.
\end{proof}

\begin{definition}[Resolvent Degree of an Algebraic Group]\label{def:RD of an Algebraic Group} Let $G$ be an algebraic group. The \textbf{resolvent degree} of $G$ over $k$ is
\begin{equation*}
	\RD_k(G) := \sup\lb \RD_k(X \bra X/G) \st X \text{ is a primitive, generically free $G$-variety over } k \rb.
\end{equation*}
\end{definition}

In the course of investigating the resolvent degree of an algebraic group, it is often useful to pass to an extension $K'/K$ of bounded resolvent degree, or more generally to an extension $K'/K$ in some specified class of accessory irrationalities $\mathcal{E}$; following \cite{ArnoldShimura1976,Reichstein2022} we can formalize the maximal such extension as follows.

\begin{definition}[$\mathcal{E}$-Closure]\label{def:RDKD closure}
Let $\mathcal{E}\colon \Fieldsk\to \Cat$ be a saturated class of accessory irrationalities. Let $K$ be a $k$-field and fix an algebraic closure $K \injects \overline{K}$. Consider the set 
\[
    \mathcal{S}_\mathcal{E}:=\{K \injects K' \injects \overline{K} \st K' \in \mathcal{E}(K)\}.
\]
We define an \textbf{$\mathcal{E}$-closure of $K$} to be the compositum
\begin{equation*}
	\KE := K\lp \mathcal{S}_\mathcal{E} \rp \injects \overline{K}.
\end{equation*}

We say that $K$ is \textbf{$\mathcal{E}$-closed} if $K = \KE$. Note that, by definition, for any finite extension of $k$-fields $K\subset L\subset \overline{K}$, if $L\in \mathcal{E}(K)$ then $L\subset\KE$. Conversely, because $\mathcal{E}$ is saturated, if $L\subset \KE$, then $L\in\mathcal{E}(K)$.

\end{definition}

\begin{example}[Common $\mce$-Closures]\label{ex:Common mce-Closures} \mbox{}
    \begin{enumerate}
        \item For $\mathcal{E}=\Ab$, $K \injects K^{\Ab}$ is the usual abelian closure (i.e. maximal abelian extension).
        \item For $\mathcal{E}=\Sol$, $K \injects \KS$ is the usual solvable closure.
        \item For $\mathcal{E}=\RD_k^{\le d}$, we write $\KD:=K^{\RD_k^{\le d}}$. This closure was first considered in \cite{ArnoldShimura1976} and studied in some detail recently in \cite{Reichstein2022}. For $d=0$, $K=K^{\lp 0\rp}$. As radicals have resolvent degree 1, we see $\KS \injects \KD$ whenever $d\ge 1$. If $d\le d'$, then $\RD_k^{\le d}$ is a subfunctor of $\RD_k^{\le d'}$, and one sees that $\KD \injects K^{\lp d'\rp}$ as expected. For more details, see \cite[Section 6]{Reichstein2022}.    
    \end{enumerate}
\end{example}

\subsubsection{Resolvent Degree of a Functor}
Reichstein extended the notion of resolvent degree to ``split'' functors in \cite[Section 8]{Reichstein2021} (see also \cite[Section 7]{Reichstein2022}). In this subsection, we extend this to a definition of resolvent degree for arbitrary functors (Definition~\ref{d:rdfun}), and we show in Lemma~\ref{l:rdgcomp} that this definition recovers Reichstein's definition if the functor is split. In particular, this establishes the equivalence of Definition~\ref{def:RD of an Algebraic Group} and Reichstein's definition of the resolvent degree of an algebraic group \cite[Definition 10.1]{Reichstein2022}.

We begin by defining the resolvent degree of a functor $F\colon \Fieldsk\to \Sets$, building on Merkurjev's definition of essential dimension in this context \cite{BerhuyFavi2003}. 
	
\begin{definition}[Merkurjev]\label{d:edfun}
	Let $F\colon \Fieldsk\to\Sets$ be a functor. Let $K$ be a $k$-field, and let $\alpha\in F(K)$. Define the {\bf essential dimension} of $\alpha$ by
	\[
	   \ed_k(\alpha):=\min \{ \trdeg_k L~|~\alpha\in \im(F(L)\to F(K))\}.
	\]
	Define the {\bf essential dimension} of $F$ by
	\[
		\ed_k(F):=\sup_{K,\alpha\in F(K)} \ed_k(\alpha).
	\]
\end{definition}

\begin{example}[Essential Dimension of $\Fin$]
    Let $F=\Fin$ as above.  Let $K$ be a $k$-field and let $\alpha\colon K\hookrightarrow L$ be a finite extension of $K$, and let $Y\dashrightarrow X$ be any branched cover of $k$-varieties such that 
    \[
    (k(X)\hookrightarrow k(Y)) \cong (K\hookrightarrow L).
    \]

    \noindent
    Tracing through the definitions, one immediately obtains that
    \[
        \ed_k(\alpha)=\ed_k(L/K)=\ed_k(Y\dashrightarrow X)
    \]
    where the left-hand side denotes the quantity defined in Definition~\ref{d:edfun}, the middle term denotes the essential dimension of a finite extension of $k$-fields as in \cite[Definition 2.1]{BuhlerReichstein1997}, and the right-hand side denotes the quantity of Definition~\ref{d:edbra}.
\end{example}

Also recall the field theoretic formulation of Definition~\ref{d:rdbra} \cite{Brauer1975}, which by \cite[Propositon 2.4]{FarbWolfson2019} is equivalent to Definition~\ref{d:rdbra} under the assignment $X\mapsto k(X)$:
\begin{definition}[Brauer]\label{d:brauer}
	Let $K\hookrightarrow L$ be a finite extension of $k$-fields.  The {\bf resolvent degree} of $L$ over $K$, $\RD_k(L/K)$ is the minimal $d$ for which there exists a finite tower of finite extensions of $k$-fields
	\[
	       K=K_0\hookrightarrow K_1\hookrightarrow\cdots\hookrightarrow K_r
	\]
	and an embedding $L\hookrightarrow K_r$ over $K$ with $\ed_k(K_{i+1}/K_i)\le d$ for all $i$.
\end{definition}

Motivated by \cite[Proposition 2.13]{FarbWolfson2019}, we combine Brauer's and Merkurjev's definitions to obtain the following.
\begin{definition}[Resolvent Degree of a Functor]\label{d:rdfun}
    Let $F\colon \Fieldsk \to \Sets$ be a functor.  Let $K$ be a $k$-field and let $\alpha\in F(K)$. We define the {\bf resolvent degree} of $\alpha$ by
	\[
        \RD_k(\alpha):=\min_{L/K~\text{finite}} \max\{\RD_k(L/K),\ed_k(\alpha|_L)\}
	\]
	Define the {\bf resolvent degree} of $F$ by
	\[
		\RD_k(F):=\sup_{K,\alpha\in F(K)} \RD_k(\alpha).
	\]
\end{definition}

\begin{example}[Resolvent Degree of $\Fin$ and $H^1(-,G)$]\mbox{}
    \begin{enumerate}
        \item Consider the functor $\Fin\colon \Fieldsk\to\Sets$.  Let $\alpha\colon K\hookrightarrow L$ be a finite extension of $k$-fields considered as an element $\alpha\in \Fin(K)$.  Then, by inspection
        \[
            \RD_k(\alpha)=\RD_k(L/K)
        \]
        where the left-hand side is as in Definition~\ref{d:rdfun} and the right-hand side is as in Brauer's Definition~\ref{d:brauer}.    
        \item For an algebraic group $G$ over $k$, consider the functor $H^1(-,G):\Fieldsk \ra \Sets$. For each $K/k$, the elements of $H^1(K,G)$ are isomorphism classes of $G$-torsors over $\spec(K)$, where the local triviality condition is with respect to the fppf topology. Given a $G$-torsor $T\to K$ with class $\alpha\in H^1(K,G)$, let $X\dashrightarrow X/G$ be an integral model.  By inspection, we see that
        \[
            \RD_k(\alpha)=\RD_k(X\dashrightarrow X/G)
        \]
        where the left-hand side is as in Definition~\ref{d:rdfun} and the right-hand side is as in Definition~\ref{d:rdgvar}.  Similarly, 
        \[
            \RD_k(H^1(-,G))=\RD_k(G),
        \]
        i.e. for algebraic groups, the resolvent degree of $H^1(-,G)$ agrees with Definition \ref{def:RD of an Algebraic Group}.
    \end{enumerate}
\end{example}
 
In \cite{Reichstein2021}, Reichstein defines a notion of resolvent degree for certain functors $F\colon \Fieldsk \to \MSETS$ taking values in {\em pointed} sets. Following the notation of \cite{Reichstein2021}, for such $F$ and $K$ a $k$-field, let $1\in F(K)$ denote the distinguished element. 

\begin{definition}[Split Functors]
    A functor $F\colon \Fieldsk\to \MSETS$ is {\bf split} if for all $K/k$ and all $\alpha\in F(K)$, there exists a finite extension $L/K$ such that $\alpha|_L=1\in F(L)$. In such a case, we say that $\alpha$ is {\bf split} by $L/K$.
\end{definition}

\begin{definition}[Reichstein]
	Let $F\colon \Fieldsk\to \MSETS$ be a split functor. Given $\alpha\in F(K)$, define the {\bf split resolvent degree} of $\alpha$ by
	\[
	   \RD^{sp}_k(\alpha):=\min\{\RD_k(L/K)~|~\alpha~\text{is split by}~L/K\}
	\]
	Define the {\bf split resolvent degree} of $F$ by
	\[
	   \RD^{sp}_k(F):=\sup_{K,\alpha\in F(K)} \RD_k^{sp}(\alpha).
	\]
\end{definition}
	
\begin{remark}[On Split Resolvent Degree]\mbox{}
    \begin{enumerate}
		\item Given a split functor $F$, we will continue to write $\RD_k(F)$ to denote the resolvent degree of $F$ considered as a functor $F\colon \Fieldsk\to \MSETS\to \Sets$ (i.e. where we forget the distinguished element).  Below we will to show that $\RD_k(F)=\RD^{sp}_k(F)$ for any such $F$.
		\item A motivating example of a split $F$ is given by $H^1(-,G)$ for $G$ an algebraic group. As above, write $\RD^{sp}_k(G):=\RD^{sp}_k(H^1(-,G))$. As Reichstein observes \cite[Conjecture 17]{Reichstein2021} (see also \cite[Conjecture 1.4]{Reichstein2022}), a folklore conjecture implicit in work of Tits is that if $G$ is a connected complex algebraic group and $K$ is a $\CC$-field, then every $G$-torsor over $K$ splits over a solvable extension $L/K$.  In particular, Tits' conjecture implies that 
		\[
			\RD^{sp}_{\CC}(G)=\RD_{\CC}(G)\le 1
		\]
		for every connected complex algebraic group $G$. Reichstein proves unconditionally \cite[Theorem 1.1]{Reichstein2022} that 
		\[
			\RD^{sp}_{\CC}(G)\le 5.
		\]
		\item The norm-residue isomorphism theorem implies that 
		\[
		      \RD^{sp}_k(H^*(-;\mu_n))=\RD_k(H^*(-;\mu_n))=1
		\]
		for every field $k$ of characteristic prime to $n$.  In particular, torsion Galois cohomology cannot detect $\RD>1$.
    \end{enumerate}
\end{remark}
			
\begin{lemma}[Equivalence of Split Resolvent Degree for Split Functors]\label{l:rdgcomp}
    Let $F\colon \Fieldsk\to\MSETS$ be a split functor. Then
	\[
	   \RD^{sp}_k(F)=\RD_k(F).
	\]
    In particular, for $G$ an algebraic group over $k$, Definition~\ref{def:RD of an Algebraic Group} agrees with \cite[Definition 10.1]{Reichstein2022}.
\end{lemma}

\begin{proof}
    For any $k$-field $K$, $\ed_k(1_K)=0$ by definition (since $1_k|_K=1_K\in F(K)$).  Therefore, the definitions immediately give that
	\[
		\RD_k(F)\le \RD_k^{sp}(F).
	\]
	We now show the opposite inequality.  Observe that it suffices to prove that for all $k$-fields $K$ and all $\alpha\in F(K)$, we have 
	\[
		\RD^{sp}_k(\alpha)\le \RD_k(\alpha).
	\]
    By \cite[Lemma 7.6(b)]{Reichstein2022}, 
    \[
		\ed_k(\alpha)\ge \RD^{sp}_k(\alpha).
	\]
	But then, for any finite extension $L/K$, we have
	\begin{align*}
		\max\{\RD_k(L/K),\ed_k(\alpha|_L)\}&\ge \max\{\RD_k(L/K),\RD^{sp}_k(\alpha|_L)\}\\
		      &=\max\{\RD_k(L/K),\min\{\RD_k(L'/L)~|~\alpha|_L~\text{is split by}~L'/L\}\}\\
			&=\min\{\max\{\RD_k(L/K),\RD_k(L'/L)\}~|~\alpha|_L~\text{is split by}~L'/L\}\\
			&\ge\min\{\RD_k(L'/K)~|~\alpha~\text{is split by}~L'/K\}\\
            &=:\RD_k^{sp}(\alpha)
	\end{align*}
	where the final inequality follows from \cite[Lemma 2.7]{FarbWolfson2019}.  Minimizing the left hand side of the above inequality over all finite extensions $L/K$, we obtain that 
	\[
		\RD_k(\alpha)\ge\RD^{sp}_k(\alpha)
	\]
	as desired.
\end{proof}


\section{Generalized Versality and Special Points}\label{sec:Generalized Versality and Special Points}
As stated at the beginning of Section \ref{sec:Background}, we fix a ground field $k$ throughout. By Definition \ref{def:RD of an Algebraic Group}, for any algebraic group $G$ and any primitive, generically free $G$-variety $X$, we have
\begin{equation*}
	\RD_k(X \bra X/G) \leq \RD_k(G). 
\end{equation*}

\noindent
It is natural to ask for which $G$-varieties $X$ we have $\RD_k(X \bra X/G) = \RD_k(G)$. We will relate this question to the notion of versality and generalizations thereof. First, we recall the definition of a versal $G$-variety. 

\begin{definition}[Section 1 of \cite{DuncanReichstein2015}]\label{def:Versal and Weakly Versal}
Let $G$ be an algebraic group defined over $k$. We say that an irreducible, generically free $G$-variety $X$ is:
\begin{itemize}
\item \textbf{weakly versal for $G$} if for every $G$-torsor $T \ra \spec(K)$, there is a $G$-equivariant $k$-morphism $T \ra X$.
\item \textbf{versal for $G$} if every $G$-invariant open subvariety of $X$ is weakly versal.
\end{itemize}  
\end{definition}

Note that $X$ being versal for $G$ is closely related to $G \racts X$ being a generic group action; the difference being that versality does not require $X/G$ to be rational (see \cite[Remark 2.8]{DuncanReichstein2015}).

In \cite[Proposition 3.7]{FarbWolfson2019}, the authors show that for a finite group $G$ and any versal $G$-variety $X$, $\RD_k(G) = \RD_k(X \bra X/G)$. While versality is a sufficient condition, it was known classically that versality is not necessary. Indeed, Klein showed in \cite{Klein1884} that $\RD_\CC(A_5) = 1$ by using the projective representation $A_5 \racts \PP_\CC^1$ (see \cite{Morrice1956} for an English translation), despite the fact that $\PP_\CC^1 \bra \PP_\CC^1/A_5$ is not $A_5$-versal.  This motivates the following generalizations of versality.

\begin{definition}[Definition 3.8 of \cite{FarbWolfson2019} and Definition 4.4 of \cite{FarbKisinWolfson2023}]\label{def:Solvably Versal and RD-Versal}\label{def:RDKD versality}
Let $G$ be an algebraic group defined over $k$ and let $X$ be an irreducible, generically free $G$-variety. Let $\mathcal{E}$ be a class of accessory irrationalities. We say that $X$ is:
\begin{itemize}
\item \textbf{weakly $\mathcal{E}$-versal for $G$} if for every $G$-torsor $T \ra \spec(K)$, there is an extension $K\hookrightarrow K'\in \mathcal{E}(K)$ and a $G$-equivariant $k$-morphism
\begin{equation*}
	T \times_{\spec(K)} \spec(K') \ra X;
\end{equation*}
\item \textbf{$\mathcal{E}$-versal for $G$} if every $G$-invariant open subvariety of $X$ is weakly $\mathcal{E}$-versal;
\item \textbf{weakly $\RD_k$-versal for $G$} if for every $G$-torsor $T \ra \spec(K)$, there is an extension $K \injects \wt{K}$ with $\RD_k\lp K \injects \wt{K} \rp \leq \RD_k(X \bra X/G)$ and a $G$-equivariant $k$-morphism 
\begin{equation*}
	T \times_{\spec(K)} \spec\lp \wt{K} \rp \ra X;
\end{equation*}
\item \textbf{$\RD_k$-versal for $G$} if every $G$-invariant open subvariety of $X$ is weakly $\RD_k$-versal. 
\end{itemize}

\end{definition}

It is immediate that for an irreducible $G$-variety $X$, we have the implications
\begin{align*}
	\text{ $X$ is versal for $G$ } &\RA \text{ $X$ is $\mathcal{E}$-versal for $G$ for any $\mathcal{E}$.} \\
    \text{ $X$ is solvably versal for $G$ } &\RA \text{ $X$ is $\RD_k^{\le 1}$-versal for $G$.} \\
     \text{ $X$ is $\RD_k^{\le d}$-versal for $G$} &\RA \text{ $X$ is $\RD_k^{\le d'}$-versal for $G$ for any $d\le d'$.} \\
    \text{ $X$ is $\RD_k^{\le d}$-versal for $G$ for $d\le \RD_k(X\bra X/G)$}&\RA \text{ $X$ is $\RD_k$-versal for $G$. }
\end{align*}

Klein showed that while $\PP_\CC^1$ is not versal for $A_5$, it is solvably versal for $A_5$. It is currently unknown if $\RD_k$-versality is strictly weaker than solvable versality (see e.g. \cite[Problem 9.3]{ChernousovGilleReichstein2006}). Nonetheless, by \cite[Proposition 3.10]{FarbWolfson2019},  $\RD_k(X \bra X/G) = \RD_k(G)$ when $X$ is $\RD_k$-versal for $G$ and $G$ is finite; note that the proof given carries over unchanged for general smooth algebraic groups $G$.

We have relatively few techniques for showing a $G$-variety is $\RD_k$-versal (especially when it is not already versal or solvably versal), and essentially none for obstructing the existence of $\RD_k$-versal varieties of a given dimension. For example, Hilbert's Sextic Conjecture \cite{Hilbert1927} --- which remains open --- predicts that there are no $\RD_\CC$-versal $A_6$-curves.

The proof of \cite[Proposition 3.10]{FarbWolfson2019} (see also \cite[Paragraph 4.1.3 and Lemma 4.9]{FarbKisinWolfson2023}), adapted to the context of Definition \ref{def:RDKD versality}, immediately yields the following:

\begin{prop}[$\RD_k(G)$ via $\RD_k^{\leq d}$-versality]\label{prop:RDK(G) via RDKD-versality}
Let $G$ be an algebraic group over $k$. Then,
\begin{equation*}
	\RD_k(G) = \min\limits_{d \geq 0}\lb \max\lb d, \dim(X) \rb \st X \text{ is a $G$-variety which is $\RD_k^{\leq d}$-versal for } G \rb.
\end{equation*}
\end{prop}

\begin{remark}[Equivalence of Definitions]\label{rem:Equivalence of Definitions}
The style of Definition \ref{def:Solvably Versal and RD-Versal} has been chosen for ease of comparison with the literature on versality, e.g. \cite{DuncanReichstein2015}. For comparison with \cite{FarbKisinWolfson2023}, and for use in what follows, we will now give an equivalent characterization in Lemma \ref{lem:RDKD-Versality via $G$-Equivariant Rational Correspondences}. 
\end{remark}

\begin{definition}[$G$-Equivariant Rational Correspondence]\label{def:G-Equivariant Rational Correspondence}
Let $G$ be an algebraic group and take $X,Y$ to be $G$-varieties. A \textbf{$G$-equivariant rational correspondence} from $Y$ to $X$ is a $G$-invariant subvariety $C \subseteq Y \times X$ such that the projection $C/G \ra Y/G$ is a generically finite, dominant morphism.
\end{definition}

\begin{lemma}[$\mathcal{E}$-Versality via $G$-Equivariant Rational Correspondences]\label{lem:RDKD-Versality via $G$-Equivariant Rational Correspondences}
Let $G$ be an algebraic group, $X$ an irreducible, generically free $G$-variety, and $\mathcal{E}$ a saturated class of accessory irrationalities. Then $X$ is:
\begin{itemize}
\item \textbf{weakly $\mathcal{E}$-versal for $G$} if for any generically free $G$-variety $Y$, there exists a $G$-equivariant rational correspondence $C \subseteq Y \times X$ such that $C/G \ra Y/G$ is in $\mathcal{E}(Y/G)$;
\item \textbf{$\mathcal{E}$-versal for $G$} if every non-empty $G$-invariant open subvariety of $X$ is weakly $\mathcal{E}$-versal for $G$.
Moreover, these conditions are equivalent to the conditions stated in \cite[Definition 4.4]{FarbKisinWolfson2023}
\end{itemize}
\end{lemma}

\begin{proof}
    We begin by proving the equivalence of the conditions of the lemma with \cite[Definition 4.4]{FarbKisinWolfson2023}.  First, observe that any $G$-equivariant rational correspondence $C\subseteq Y\times X$ with $C/G\ra Y/G$ in $\mathcal{E}(Y/G)$ gives the data of \cite[Definition 4.4]{FarbKisinWolfson2023}; indeed, item (1) of \cite[Definition 4.4]{FarbKisinWolfson2023} is the given accessory irrationality, item (2) is the map $C/G\ra X/G$ (from projecting onto the second factor), and the isomorphism (3) $C/G\times_{X/G}X\cong C$ follows from generic freeness. 
    
    Conversely, given the data specified in \cite[Definition 4.4]{FarbKisinWolfson2023}, i.e. $E\to Y/G$ in $\mathcal{E}(Y/G)$ with $f\colon E\ra  X/G$ and $f^* X\cong E\times_{Y/G}Y$, let
    \begin{align*}
        C/G &:= \operatorname{Im}(E\to Y/G\times X/G),\\
        C &= C/G\times_{X/G} X\subset Y\times X.
    \end{align*}

    \noindent
    Then, $C$ is a $G$-invariant rational correspondence; $C/G\to Y/G$ is in $\mathcal{E}(Y/G)$, as $E\to Y$ is; and $\mathcal{E}$ is saturated. We conclude that the conditions of the lemma are in fact equivalent to those of \cite[Definition 4.4]{FarbKisinWolfson2023}.  
    
    Since \cite[Definition 4.4]{FarbKisinWolfson2023} does not mention torsors, for the sake of completeness we now show that the conditions of \cite[Definition 4.4]{FarbKisinWolfson2023} are equivalent to those of Definition~\ref{def:Solvably Versal and RD-Versal}. For this, it suffices to prove the statement about weak $\mathcal{E}$-versality (as the statement about $\mathcal{E}$-versality amounts to verifying that every dense Zariski open $U\subset X$ is weakly $\mathcal{E}$-versal).  To go from a $G$-torsor $T\to\spec(K)$ as in Definition~\ref{def:Solvably Versal and RD-Versal} to the data of \cite[Definition 4.4]{FarbKisinWolfson2023} pick an integral model; to go the other way, restrict to a generic point of an irreducible component of $Y/G$.\footnote{\cite{FarbKisinWolfson2023} writes $\tilde{Y}\to Y$ for our $Y\to Y/G$.} From the proof of Lemma~\ref{l:aifin=aibra}, we see that under this correspondence, an extension $K\hookrightarrow K'$ is in $\mathcal{E}(K)$ if and only if any integral model $E\to Y/G$ is in $\mathcal{E}(Y/G)$.  The equivalence of the two definitions now follows by inspection.
\end{proof}

In \cite{DuncanReichstein2015}, Duncan and Reichstein connect versality to existence of rational points. Specifically:

\begin{theorem}[Versality and Rational Points, Theorem 1.1 (a,b) of \cite{DuncanReichstein2015}]\label{thm:Duncan-Reichstein}
Let $G$ be a linear algebraic group over $k$ and take $X$ to be an irreducible, generically free, quasiprojective $G$-variety. Then, $X$ is:
\begin{enumerate}
\item weakly versal if and only if for every $G$-torsor $T \ra \spec(K)$, $\TX(K) \not= \emptyset$;
\item versal if and only if for every $G$-torsor $T \ra \spec(K)$, $K$-points are dense in $\TX$.
\end{enumerate}
\end{theorem}

\begin{remark}[Linear Algebraic Groups vs. Smooth Algebraic Groups]
    Duncan and Reichstein restrict their attention to {\em linear} algebraic groups $G$ (they also leave the assumption of ``generically free'' implicit in their statement). In \cite[Remark 2.6]{DuncanReichstein2015}, they state that this is ``vitally important'' for their reformulations of versality.  On the other hand, a careful reading of \cite[Sections 1-4]{DuncanReichstein2015} shows that linearity can be weakened to smoothness at the cost only of rendering their Theorem 1.1(b) (that versality is equivalent to the density of $K$-points in all twisted forms) potentially vacuous, as versal varieties do not exist for general smooth $G$. 
\end{remark}

Motivated by this result, along with \cite[Lemma 4.15]{Reichstein2022}, we establish analogous claims about $\mathcal{E}$-versality.

\begin{theorem}[Generalized Versality and Special Points]\label{thm:RDKD-versality and Special Points}
Let $G$ be an algebraic group over $k$ and let $X$ be an irreducible, generically free, quasiprojective $G$-variety. Let $\mathcal{E}$ be a saturated class of accessory irrationalities.  Then: 
\begin{enumerate}
\item $X$ is weakly $\mathcal{E}$-versal if and only if for every $G$-torsor $T\ra \spec(K)$ with $K$ finitely generated over $k$, $\TX\lp \KE\rp\not= \emptyset$.
\item If $G$ is smooth, $X$ is $\mathcal{E}$-versal if and only if for every $G$-torsor $T \ra \spec(K)$ with $K$ finitely generated over $k$, $\KE$-points are dense in $\TX$. 
\end{enumerate}
\end{theorem}

\begin{remark}[Context for Smoothness]
    Cartier showed that every algebraic group over a field $k$ of characteristic 0 is smooth (see \cite[Corollary 5.2.18]{Poonen2017}), so in this case, no assumption on $G$ is needed in the second part of the theorem.  We assume smoothness in order to reduce the proof to a Galois descent argument.  More generally, it seems reasonable to expect that the theorem holds without any assumptions on $G$ and $k$, at the cost of using fppf descent in lieu of Galois descent.  We do not pursue this here.
\end{remark}

\begin{proof}[Proof of Theorem~\ref{thm:RDKD-versality and Special Points}]
We begin by showing the first statement. Suppose that $X$ is weakly $\mathcal{E}$-versal. Let $K$ be a field which is finitely generated over $k$ and consider a $G$-torsor $T \ra \spec(K)$ with integral model $Y \ra Y/G$. The $G$-equivariant isomorphism $Y \times_{Y/G} \spec(K) \cong T$ induces
\begin{equation*}
	(Y \times X)/G \times_{Y/G} \spec(K) \cong \TX.
\end{equation*}

\noindent
As $X$ is weakly $\mathcal{E}$-versal for $G$, there exists a $G$-equivariant rational correspondence $C \subseteq Y \times X$ with $C/G \bra Y/G$ in $\mathcal{E}(Y/G)$. Taking the quotient of the inclusion $C \ra Y \times X$ yields the morphism $C/G \ra (Y \times X)/G$ of $Y/G$-varieties. We can restrict the morphisms $C/G \ra Y/G$ and $(Y \times X)/G \ra Y/G$ along the generic point $\spec(K) \ra Y/G$ and obtain the following morphism of pullback diagrams, where $\spec(L)$ denotes the generic point of $C/G$:

\begin{equation*}
    \xymatrix{
        \spec(L) \ar[rr] \ar[dr] \ar[dd] && \spec(K) \ar@{-}[d] \ar@{=}[dr] \\
        & \TX \ar[rr] \ar[dd] & \ar[d] & \spec(K) \ar[dd] \\
        C/G \ar@{-}[r] \ar[dr] & \ar[r]  & Y/G \ar@{=}[dr] \\
        & (Y\times X)/G \ar[rr] && Y/G
    }
\end{equation*}

Since $C/G \ra Y/G$ is in $\mathcal{E}(Y/G)$ (by assumption), $L\in\mathcal{E}(K)$ (by \ref{it:accirrbirat} of Definition~\ref{d:aibra}). From Definition \ref{def:RDKD closure}, there is an inclusion of $K$-fields $L \injects \KE$. Consequently, the constructed $\KE$-point
\begin{equation*}
	\spec\lp \KE \rp \ra \spec(L) \ra \TX,
\end{equation*}

\noindent
shows that $\TX\lp \KE \rp \not= \emptyset$. 

Now, suppose that $T \ra \spec(K)$ and $Y \ra Y/G$ are as above and $\spec\lp \KE \rp \ra \TX$ is a $\KE$-point. Since $\TX$ is a $K$-variety, there is a finite extension $K \injects E$ such that $\spec\lp \KE \rp \ra \TX$ factors as
\begin{equation*}
	\spec\lp \KE \rp \ra \spec(E) \ra \TX.
\end{equation*}
As $\mathcal{E}$ is saturated, $E\in\mathcal{E}(K)$.
\noindent
Now, let $C/G$ denote the closure of $\spec(E)$ in $(Y \times X)/G$ and take $C \subseteq Y \times X$ to be the preimage of $C/G$ under the quotient map $Y \times X \ra (Y \times X)/G$. Consequently, $C$ is a $G$-invariant subvariety, $C \ra Y$ is a $G$-equivariant, generically finite, dominant morphism (by construction), and $C/G \bra Y/G$ is in $\mathcal{E}(Y/G)$ by construction. From Lemma \ref{lem:RDKD-Versality via $G$-Equivariant Rational Correspondences}, we see that $X$ is weakly $\mathcal{E}$-versal. 

It remains to show the second claim.  First, note that for any variety $Z$ and field $\wt{K}$, $\wt{K}$-points are dense in $Z$ if and only if each Zariski open of $Z$ contains a $\wt{K}$-point. Next, in the setting of the theorem, for a Zariski open $V\subset \TX$, consider the quotient map $q\colon T\times X\to \TX$, along with the projection map $p_X\colon T\times X\to X$, and set $U=p_X(q^{-1}(V))\subset X$. 

By inspection, $U\subset X$ is a $G$-invariant Zariski open, and $V\subset \prescript{T}{}{U}$.  Therefore, if $\TX\lp\KE\rp$ is dense in $\TX$, then $\prescript{T}{}{U}\lp\KE\rp\neq \emptyset$ for every $G$-invariant Zariski open $U\subset X$. From our argument above, we conclude that every such $U$ is weakly $\mathcal{E}$-versal for $G$, and thus $X$ is $\mathcal{E}$-versal for $G$.  It remains to show the converse.

Suppose that $X$ is $\mathcal{E}$-versal for $G$. We need to show that $\TX\lp \KE\rp$ is dense.  As $G$ is smooth, we can identify \'etale and fppf cohomology with coefficients in $G$ (see e.g. \cite[Remark III.4.8(a)]{Milne1980}); this allows us to make arguments via Galois descent. 

Recall that a (right) $G$-torsor $T\to \spec(K)$ with cocycle $\tau\in H^1(K,G)$ is also a (left) $G^\tau$-torsor, where $G^\tau$ is the inner twist of $G$ over $K$ determined by $\tau$. We can see this explicitly as follows.  Indeed, let $\KSEP$ denote a separable closure of $K$. Then, letting $\Gal(K):=\Aut_K(\KSEP)$, we obtain a 1-cocycle 
\[
    \tau\colon \Gal(K)\to G(\KSEP)
\]
from $T$ by picking an element $0\in T(\KSEP)$.  Explicitly, $\tau$ is defined to be the map such that for $\sigma\in\Gal(K)$
\[
    \prescript{\sigma}{}{0}=0\cdot \tau(\sigma)
\]
where we write $\prescript{\sigma}{}{x}$ to denote the $\sigma$-translate of a $\KSEP$ point $x$ of a $K$-variety.  Note that $\tau(\sigma)$ is uniquely determined because $T(\KSEP)$ is a principal right $G(\KSEP)$-set.  Further, the choice of $0$ determines an isomorphism of right $G(\KSEP)$-sets
\begin{align*}
    \varphi_0\colon G(\KSEP)&\stackrel{\cong}{\rightarrow} T(\KSEP)\\
    g & \mapsto 0\cdot g
\end{align*}
The left action of $G(\KSEP)$ on itself now defines a left $G(\KSEP)$-action on $T(\KSEP)$ via
\[
    g\cdot (0\cdot h):=0\cdot gh.
\]
By inspection, this is not equivariant for the standard $\Gal(K)$-action on $G(\KSEP)$, but rather for the $\tau$-twisted action
\[
    \sigma\cdot g:= \tau(\sigma) \lp \prescript{\sigma}{}{g} \rp\lp \tau(\sigma)^{-1} \rp.
\]
By Galois descent, just as in \cite[5.12.5.1]{Poonen2017}, we conclude that $T\to \spec(K)$ is actually a $G^\tau-G$ bitorsor, and therefore, the twist $\TX$ carries a left action of $G^\tau$.  

Now, given a Zariski open $V\subset \TX$ as above, consider the Zariski open $G^\tau \cdot V\subset \TX$.  Then $V\lp \KE\rp\neq \emptyset$ if and only if $G^\tau \cdot V\lp \KE\rp\neq\emptyset$. Now let $q\colon T\times X\to \TX$ denote the quotient map, let $p_X\colon T\times X\to X$ denote the projection, and set $U := p_X(q^{-1}(V)) \subseteq X$.  Note that $U$ is a $G$-invariant Zariski open by inspection.  Since $X$ is $\mathcal{E}$-versal for $G$, $U$ is weakly $\mathcal{E}$-versal.  By the proof of the first statement of the theorem above, we have that $\prescript{T}{}{U}\lp \KE\rp\neq \emptyset$. To conclude, we claim that 
\[
    G^\tau\cdot V=\prescript{T}{}{U}.
\]
Granting this claim, we have, by the above, that $V\lp\KE\rp\neq\emptyset$, and thus that $\TX\lp\KE\rp$ is dense in $\TX$, as claimed.  We prove the claim by a straightforward Galois descent argument.  Indeed, our choice of $0\in T(\KSEP)$ determines an isomorphism
\begin{align*}
    \prescript{T}{}{U}(\KSEP)&\cong\{(0\cdot g,u)\in G(\KSEP)\times U(\KSEP)\}/(0\cdot g,u)\sim (0,g\cdot u)\\
    &=\{(0\cdot g,x)\in G(\KSEP)\times X(\KSEP)~|~\exists [(0\cdot h,v)]\in V(\KSEP) \text{ s.t. } h\cdot v=x\}/(0\cdot g,x)\sim (0,g\cdot x)\\
    &=\{(0\cdot gh,v)\in G(\KSEP)\times X(\KSEP)~|~[(0\cdot h,v)]\in V(\KSEP)\}/(0\cdot g,x)\sim (0,g\cdot x)\\
    &=(G^\tau\cdot V)(\KSEP).
\end{align*}
By inspection, this isomorphism is $\Gal(K)$-equivariant, and thus $\prescript{T}{}{U}=G^\tau\cdot V$ as claimed.
\end{proof}

Taken together, Proposition \ref{prop:RDK(G) via RDKD-versality} and Theorem \ref{thm:RDKD-versality and Special Points} allow us to re-contextualize problems about the resolvent degree of algebraic groups as questions about special points on twists of $G$-varieties. Indeed, we can ask about necessary and sufficient conditions for a variety over $K$ to have an $\mathcal{E}$-point or a dense collection of $\mathcal{E}$-points. As an example, we have the following, which appeared as \cite[Lemma 4.7]{FarbKisinWolfson2023} for the case of finite $G$.

\begin{lemma}[$\mathcal{E}$-Versality and Composition]\label{lem:RDKD-versality and Composition}
Let $G$ be a  smooth algebraic group over $k$, let $\mathcal{E}$ be a saturated class of accessory irrationalities which is closed under extensions, and let $Y$ be a primitive $G$-variety which is $\mathcal{E}$-versal for $G$. Suppose $X$ is an irreducible, generically free $G$-variety which admits a $G$-equivariant rational correspondence $C\subset Y\times X$ such that $C\ra X$ is dominant and $C/G\ra Y/G$ is in $\mathcal{E}(Y/G)$. Then $X$ is $\mathcal{E}$-versal.
\end{lemma}

\begin{remark}[Representatives for $\mathcal{E}$-Versality]
    In contrast to the definition of $\mathcal{E}$-versality, Lemma \ref{lem:RDKD-versality and Composition} allows us to test $\mathcal{E}$-versality of $X$ by looking at correspondences from the {\em single} $G$-variety $Y$, rather than from {\em all} $G$-varieties.
\end{remark}

\begin{proof}[Proof of Lemma \ref{lem:RDKD-versality and Composition}]
For weak $\mathcal{E}$-versality, let $Z$ be a generically free $G$-variety. From Lemma \ref{lem:RDKD-Versality via $G$-Equivariant Rational Correspondences}, there is a $G$-equivariant rational correspondence $D \subseteq Z \times Y$ with $D/G \ra Z/G$ in $\mathcal{E}(Z/G)$. Then, we can consider $C \times_Y D \subseteq Z \times X$, which is $G$-invariant, and we see that $(C\times_Y D)/G\ra D/G\ra Z/G$ is in $\mathcal{E}(Z/G)$ because $(C\times_YD)/G\ra D/G$ is in $\mathcal{E}(D/G)$, $D/G\ra Z/G$ is in $\mathcal{E}(Z/G)$, and $\mathcal{E}$ is closed under extensions.

For $\mathcal{E}$-versality, by Theorem~\ref{thm:RDKD-versality and Special Points}, it suffices to show that $\TX(\KE)$ is dense in $\TX$ for any $G$-torsor $T\to \spec(K)$. Fix such a $T$.  By Theorem~\ref{thm:RDKD-versality and Special Points}, $\prescript{T}{}{Y}(\KE)$ is dense in $\prescript{T}{}{Y}$. By assumption, $C\to Y$ is generically finite, dominant and $C/G\to Y/G$ is in $\mathcal{E}(Y/G)$. Therefore, by the definition of $\KE$ as an $\mathcal{E}$-closure, the density of $\prescript{T}{}{Y}(\KE)$ in $\prescript{T}{}{Y}$ implies that $\prescript{T}{}{C}\lp \KE\rp$ is dense in $\prescript{T}{}{C}$ as well. But, the map $\prescript{T}{}{C}\to \TX$ is dominant, by assumption, so we conclude that $\TX\lp\KE\rp$ is dense in $\TX$ as claimed.
\end{proof}

Theorem \ref{thm:RDKD-versality and Special Points} connects (weak) $\mathcal{E}$-versality of $X$ for $G$ to existence of special points (i.e. $\KE$-points), but still requires one to consider twists of $X$ by all torsors over finitely generated $k$-fields. Just as Proposition \ref{prop:RDK(G) via RDKD-versality} allows us to reduce from all $G$-varieties $X$ to those which are $\mathcal{E}$-versal, the following definition will give us the language to reduce the class of torsors one must consider.

\begin{definition}[$\mathcal{E}$-Versality for $G$-torsors]\label{def:RDKD-versality for G-torsors}
Let $G$ be an algebraic group over $k$. Let $\mathcal{E}$ be a class of accessory irrationalities. A $G$-torsor $T \ra \spec(K)$ is \textbf{$\mathcal{E}$-versal} for $G$ if there exists any integral model $Y \ra Y/G$ of $T \ra \spec(K)$ such that $Y$ is an $\mathcal{E}$-versal $G$-variety. 
\end{definition}

Note that, because all integral models are birational to each other, if $T\ra \spec(K)$ is $\mathcal{E}$-versal, then every integral model $Y\ra Y/G$ is $\mathcal{E}$-versal. We can now restate a variant of Lemma \ref{lem:RDKD-versality and Composition} as follows:

\begin{lemma}[An Equivalent Version of Lemma \ref{lem:RDKD-versality and Composition}]\label{lem:An Equivalent Version of RDKD-versality and Composition}
Let $G$ be a smooth algebraic group over $k$. Let $\mathcal{E}$ be a saturated class of accessory irrationalities that is closed under extensions. Suppose that $T \ra \spec(K)$ is a $G$-torsor which is $\mathcal{E}$-versal for $G$ and $X$ is a generically free $G$-variety. Then, $X$ is $\mathcal{E}$-versal if and only if $\TX$ has a dense collection of $\KE$-points.
\end{lemma}

\begin{proof}
From Theorem \ref{thm:RDKD-versality and Special Points}, if $X$ is $\mathcal{E}$-versal, then $\TX$ has a dense collection of $\KE$-points.

Now, suppose that $\TX$ has a dense collection of $\KE$-points and let $Y \ra Y/G$ be an integral model of $T \ra \spec(K)$ (which, as remarked above, is $\mathcal{E}$-versal because $T\ra\spec(K)$ is). Let $\spec\lp \KE \rp \ra \TX$ be a $\KE$-point of $\TX$. As established in the proof of Theorem \ref{thm:RDKD-versality and Special Points}, $\TX$ is the generic fiber of $(Y \times X)/G \ra Y/G$ and we take $C/G \subseteq (Y \times X)/G$ to be the closure of our $
\KE$-point $\spec\lp \KE \rp \ra \TX$. We set $C = C/G \times_{(Y \times X)/G} (Y \times X)$ and observe that $C \subseteq Y \times X$ is a $G$-equivariant rational correspondence with $C/G \ra Y/G$ in $\mathcal{E}(Y/G)$, by construction. By assumption, $Y \ra Y/G$ is $\mathcal{E}$-versal and thus $X$ is $\mathcal{E}$-versal by Lemma \ref{lem:RDKD-versality and Composition}.
\end{proof}

In another direction, we also have:

\begin{lemma}[$\mathcal{E}$-Versality for Non-Abelian, Finite, Simple Groups]\label{lem:RDKD-versality for non-abelian, finite, simple groups}
Assume $\text{char}(k)=0$. Let $G$ be a non-abelian, finite, simple group $G$, let $\mathcal{E}$ be a class of accessory irrationalities, and let $X$ be a smooth, irreducible, generically free $G$-curve over $k$. Then $X$ is weakly $\mathcal{E}$-versal for $G$ if and only if $X$ is $\mathcal{E}$-versal for $G$.
\end{lemma}

\begin{proof}
By definition, $\mathcal{E}$-versality immediately implies weak $\mathcal{E}$-versality. Now, suppose that $X$ is weakly $\mathcal{E}$-versal. Let $Y$ be a generically free $G$-variety with $\dim(Y) \geq 1$ and consider a $G$-equivariant rational correspondence $C \subseteq Y \times X$ with $C/G \ra Y/G$ in $\mathcal{E}(Y/G)$. It suffices to show that the $G$-equivariant map $C \ra X$ is dominant. Denote the scheme-theoretic image of $C$ in $X$ by $Z$. Then, $Z$ is an irreducible $G$-invariant subscheme. If $\dim(Z)=1$, then the map is dominant and we are done. Now, suppose that $\dim(Z)=0$. Then, $Z^\text{red} \in X$ is a fixed point for $G$. However, the stabilizer of any point in $X$ is abelian \cite[Theorem 1.1]{ReichsteinYoussin2000}, hence there are no fixed points for $G$ and thus $\dim(Z) \geq 1$. 
\end{proof}

Let us now revisit Hilbert's Sextic Conjecture. In light of Proposition \ref{prop:RDK(G) via RDKD-versality}, Theorem \ref{thm:RDKD-versality and Special Points}, and Lemma \ref{lem:RDKD-versality for non-abelian, finite, simple groups}, we can re-state the conjecture as follows:

\begin{conjecture}[Hilbert's Sextic Conjecture]\label{conj:Hilbert's Sextic Conjecture}
Let $T \ra \spec\lp \CC(x,y) \rp$ be the $A_6$-torsor associated to the Valentiner action $A_6 \racts \PP_\CC^2$. For any smooth, irreducible, generically free $A_6$-curve $X$, $\TX\lp \CC(x,y)^{(1)} \rp = \emptyset$.
\end{conjecture}

\noindent
Similarly, \cite[Problem 9.3]{ChernousovGilleReichstein2006} for $A_6$ can be re-stated as:

\begin{problem}[Chernousov-Gille-Reichstein]\label{prob:Chernousov-Gille-Reichstein}
Does there exist a smooth, irreducible, generically free $A_6$-curve $X$ with $\TX\lp \CC(x,y)^{\Sol} \rp \not= \emptyset$?
\end{problem}

More generally, showing that $\RD_k(G)>1$ is a question of \emph{obstructing} the existence of $K^{(1)}$-points for a sufficient supply of $K$-curves (namely, twists of $G$-curves over $k$). Many tools in the literature for obstructing rational points appear to be inadequate for this. For example, by \cite[Theorem 16.1]{MerkurjevSuslin1983}, the Brauer group of a solvably closed field is trivial.  More generally, the same holds for $H^1(-,G)$ for any connected algebraic group $G$ without simple factors of type $E_8$ (see \cite[Theorem 1.1]{Reichstein2022}); conjecturally, the same holds for all connected algebraic groups (see \cite[Conjecture 1.4]{Reichstein2022} and the discussion just preceding it). It would be instructive to turn this observation (that $G$ connected implies $H^1(\KS,G)=1$) into a proof that Brauer-Manin style invariants are insufficient to obstruct solvable points on varieties over e.g. two-dimensional $\CC$-fields. We echo Poonen's view in \cite[p. 257]{Poonen2017} that ``We need some new obstructions!''


\section{Resolvent Degree and the Sporadic Groups}\label{sec:Resolvent Degree and the Sporadic Groups}

\subsection{Upper Bounds on the Resolvent Degree of the Sporadic Groups}\label{subsec:Upper Bounds on the Resolvent Degree of the Sporadic Groups}

Recall that the Classification of Finite Simple Groups consists of 18 infinite families and 26 sporadic groups. The 26 sporadic groups are often organized as in Figure \ref{table:Historical Organization of Sporadic Groups}.

\begin{figure}[ht!]
\begin{center}
\caption{Historical Organization of Sporadic Groups}
\label{table:Historical Organization of Sporadic Groups}
\begin{tabular}{|c|c|c|c|}
\hline
Cluster &Generation &Description &Groups\\
\hline
Happy Family &First   &The Mathieu Groups &$\MOO, \MOW, \MWW, \MWH, \MWF$,\\
\hline
Happy Family &Second  &The Leech Lattice Groups &$\COO, \COW, \COH, \SUZ, \MCL, \HS, \JW$,\\
\hline
Happy Family &Third   &Other Monster Subgroups &$\FIWW, \FIWH, \FIWF, \THO, \HN, \HE, \B, \MO$,\\
\hline
The Pariahs  &        &The Pariahs &$\JO, \JH, \JF, \ON, \RU, \LY$.\\
\hline
\end{tabular}
\end{center}
\end{figure}

\noindent
From \cite[Theorems 1.2 and 1.3]{Reichstein2022}, $\RD_k(G) \leq \RDC(G)$ for any finite simple group. Thus, to determine upper bounds on the resolvent degree of the sporadic groups, it suffices to work over $\CC$. For each sporadic group $G$, we will determine a complex $G$-variety $X_G$ such that $\RDC(G) \leq \dim_{\CC}\lp X_G \rp$. We begin by considering a minimal dimensional projective representation. It is immediate that any linear representation of $G$ yields a projective representation of $G$, however these are not the only projective representations of $G$. Indeed, there are groups $\Gamma$ such that the projectivizations of linear representations of $\Gamma$ correspond exactly to projective representations of $G$. Such a group $\Gamma$ is called a \textbf{Schur cover} of $G$ (or sometimes a \textbf{Schur representation group} of $G$).  Each sporadic group $G$ is perfect, hence the Schur covers of $G$ are isomorphic and so we simply refer to \emph{the} Schur cover of $G$ henceforth. Explicitly, the Schur cover of $G$ is a central extension of $G$ by the \textbf{Schur multiplier} $\SCH(G) = H^2\lp G, \CC^* \rp$, which is a finite abelian group whose exponent divides the order of $G$. For more on projective representations of finite groups, we refer the reader to \cite[Chapter 11]{Isaacs1976}. 

Given $G$ and a projective representation $\PP(\rho):G \ra \PGL(V)$ coming from a linear representation $\rho$ of the Schur cover, we are not interested in just $\PP\lp V \rp$, but $G$-invariant subvarieties thereof. We can construct such invariant subvarieties by looking at the vanishing of $G$-invariant polynomials. Note that the vector space of homogeneous polynomials of degree $d$ which are invariant under the Schur cover of $G$ is $\Sym_{a.G}^d\lp V^\vee \rp$. We set $m_d(\rho) = \dim\lp \Sym_{a.G}^d\lp V^\vee \rp \rp$ and note that the \textbf{Molien series} of $\rho$ is the generating function
\begin{equation*}
M(\rho;t) := \SL_{d \geq 0} m_d(\rho)t^d.
\end{equation*}

We now introduce notation for the minimal projective representations for each sporadic group and address computing the relevant Molien series.

\begin{notation}[Representations of Sporadic Groups]\label{not:Representations of Sporadic Groups}
For the groups $G=\MOO,$ $\MWH,$ $\MWF,$ $\JO,$ $\JF,$ $\COH,$ $\COW,$ $\FIWW,$ $\FIWH,$ $\HS,$ $\MCL,$ $\HE,$ $\HN,$ $\THO,$ $\LY,$ $\B$, and $\MO$, a projective representation of minimal dimension arises as the projectivization of an irreducible linear representation of $G$. When $G$ is not one of $\HS$, $\MCL$, $\FIWW$, or $\B$, this claim is immediate, as the Schur multiplier is trivial. For the cases $G=\HS$, $\MCL$, $\FIWW$, and $\B$, one can verify this claim by inspecting the character tables of the Schur covers. We set $d(G)$ to be the minimal dimension of a non-trivial linear representation of $G$ and $\rho_G$ to be the representation corresponding to the first character $\chi$ in the $\ATLAS$ character table for $G$ for which $\chi(1)=d(G)$ (note that $\rho_{G}$ is necessarily irreducible by our minimality assumption). Additionally, set $V_G$ to be the vector space corresponding to $\rho_G$.

For the groups $G=\MWW$, $\JW$, $\JH$, $\COO$, $\FIWF$, $\SUZ$, $\RU$, and $\ON$, a projective representation of minimal dimension only arises as the projectivization of an irreducible linear representation of the Schur cover of $G$. Correspondingly, we set $a(G) = \lav \SCH(G) \rav$, $d(G)$ to be the minimal dimension of a non-trivial representation of the Schur cover (denoted by $a.G$), $\rho_{G}$ to be the representation corresponding to the first character in the $\ATLAS$ character table for $a.G$ for which $\chi(1) = d(G)$ (as above, our minimality assumption guarantees that $\rho_G$ is irreducible), and $V_{G}$ to be the vector space corresponding to $\rho_{G}$.

When $G$ is clear from the context, we simply write $a$ and $d$. Additionally, we note that the order of the characters in the $\ATLAS$ \cite{ConwayCurtisNortonParkerWilson1985} is the same as in \GAP \ character table library \cite{Breuer2013}. 
\end{notation}

\begin{remark}[Schur Multiplier $\MWW$]\label{rem:Schur Multiplier of MWW}
    In \cite[p.739-741]{BurgoyneFong1966}, it is incorrectly claimed that $a\lp \MWW \rp = 3$. In the correction \cite{BurgoyneFong1968},  it is incorrectly asserted that $a\lp \MWW \rp = 6$. Finally, \cite[Section V]{Mazet1979} correctly establishes that $\SCH\lp \MWW \rp \cong \ZZ/12\ZZ$.
\end{remark}

\begin{remark}[The Unique Case of $\MOW$]\label{rem:The Unique Case of MOW}
    There is another error in \cite[p.739-741]{BurgoyneFong1966}, where they incorrectly claim that $a\lp \MOW \rp=1$. However, \cite{BurgoyneFong1968} correctly establishes that $a\lp \MOW \rp = 2$.

    An observant reader may have noticed that $\MOW$ does not appear in Notation \ref{not:Representations of Sporadic Groups}. While we will use the same notation conventions, $\rho_{\MOW}$ is not a projective representation of minimal dimension for $\MOW$. While the Schur cover $2.\MOW$ admits a 10-dimensional linear representation, we will instead take $\rho_{\MOW}$ to be the first 11-dimensional linear representation of $\rho_{\MOW}$. We will justify this choice in Remark \ref{rem:The Unique Case of MOW, II}, after we prove Theorem \ref{thm:Bounds on the Resolvent Degree of the Sporadic Groups}.
\end{remark}

\begin{remark}[Computation of Molien Series / Molien Series Coefficients]
For the sporadic groups $G$ where $|G|$ and $d(G)$ are sufficiently small ($\MOO,$ $\MOW,$ $\MWW,$ $\MWH,$ $\MWF,$ $\JO,$ $\JW$, $\JH$, $\COH,$ $\COW,$ $\COO$, $\SUZ$, $\HS,$ $\MCL,$ $\RU$, and $\HE$), we compute the Molien series $M\lp \rho_G;t \rp$ as a rational function using the character table library \cite{Breuer2013} in \SAGE \ by accessing \GAP. For the remaining groups $G$ ($\JF,$ $\FIWW,$ $\FIWH,$ $\FIWF,$ $\HN,$ $\THO,$ $\ON$, $\LY,$ $\B$, and $\MO$), we store data for $\rho_G$ from the character table library \cite{Breuer2013}, which we then use to compute the first 20 coefficients $m_1\lp \rho_G \rp, \dotsc, m_{20}\lp \rho_G \rp$ of $M\lp \rho_G;t \rp$ in $\SAGE$. The \SAGE \ script, data files (for $\JF,$ $\FIWW,$ $\FIWH,$ $\FIWF,$ $\HN,$ $\THO,$ $\ON$, $\LY,$ $\B$, and $\MO$), and output files (for all sporadic groups), are available at
\begin{center}
\href{https://www.alexandersutherland.com/research/RD-for-the-Sporadic-Groups}{\texttt{https://www.alexandersutherland.com/research/RD-for-the-Sporadic-Groups}},
\end{center}

\noindent
or by downloading the source package on the \arxiv \ version of this work. Additionally, for every sporadic group $G$, we record the beginning of the power series expansion of $M\lp \rho_G;t \rp$ in Appendix \ref{app:Power Series Expanions of Molien Series}.
\end{remark}

For each sporadic group $G$, Figure \ref{table:Dimensions of Projective Representations and Degrees of Invariants} records $\dim\lp \PP\lp V_G \rp \rp$ and a list of the degrees of the invariants we will use in what follows. Note that ordering of the groups in Figure \ref{table:Dimensions of Projective Representations and Degrees of Invariants} is determined by $\dim\lp \PP\lp V_G \rp \rp$. Our construction for $\JW$ does not require any invariants, so we leave the corresponding entry blank. Additionally, our proof will use the minimal dimensional permutation representation for the sporadic groups other than $\MOO$, $\MOW$, $\MWH$, and $\MWF$, hence Figure \ref{table:Dimensions of Projective Representations and Degrees of Invariants} includes $\dim\lp \PERM_G \rp$ as well. For each $G$, we denote this representation by $\PERM_G$ and observe that $\dim\lp \PERM_G \rp = [G:H]$, where $H$ is a maximal subgroup of $G$ of maximal order. For details on maximal subgroups, see \cite{DietrichLeePopiel2023} when $G=\MO$, \cite{Wilson2009} when $G=\FIWW$, $\FIWH$, $\FIWF$, $\JF$, $\THO$, and $\B$, \cite{Wilson2017} for a survey, and \cite{ConwayCurtisNortonParkerWilson1985}) for all other cases.

\newpage

\begin{figure}[ht!]
\begin{center}
\caption{Dimensions of Projective Representations and Degrees of Invariants}
\label{table:Dimensions of Projective Representations and Degrees of Invariants}
\begin{tabular}{|c|c|c|c|}
\hline
Group $G$ &$\dim\lp \PP\lp V_G \rp \rp$ &Degrees of Relevant Invariants &$\dim\lp \PERM_G \rp$\\
\hline
$\JW$ &5 &N/A &100\\
\hline
$\MOO$ &9 &2, 3, 4 &11\\
\hline
$\MWW$ &9 &4 &22\\
\hline
$\MOW$ &10 &2, 3, 4 &12 \\
\hline
$\SUZ$ &11 &12 &1782\\
\hline
$\JH$ &17 &6 &85\\
\hline
$\MWH$ &21 &2, 3, 4, 5 &23\\
\hline
$\HS$ &21 &2, 4, 5 &100\\
\hline
$\MCL$ &21 &2, 5 &275\\
\hline
$\MWF$ &22 &2, 3, 4, 5 &24\\
\hline
$\COH$ &22 &2, 6 &276\\
\hline
$\COW$ &22 &2, 8 &2300\\
\hline
$\COO$ &23 &2, 12 &98280\\
\hline
$\RU$ &27 &4 &4060\\
\hline
$\HE$ &50 &3, 4 &2058\\
\hline
$\JO$ &55 &2, 3, 4, 4 &266\\
\hline
$\FIWW$ &77 &2, 6, 8 &3510\\
\hline
$\HN$ &132 &2, 6, 7 &1140000\\
\hline
$\THO$ &247 &2, 8, 8 &143127000\\
\hline
$\ON$ &341 &6, 6, 6 &122760\\
\hline
$\FIWH$ &781 &2, 3, 4, 5, 5 &31671\\
\hline
$\FIWF$ &782 &3, 6, 6 &306936\\
\hline
$\JF$ &1332 &4, 6, 6, 7 &173067389\\
\hline
$\LY$ &2479 &6, 6, 6, 6 &8835156\\
\hline
$\B$ &4370 &2, 4, 6, 8, 8 &13571955000\\
\hline
$\MO$ &196882 &2, 3, 4, 5, 6, 6, 6, 7 &97239461142009186000\\
\hline
\end{tabular}
\end{center}
\end{figure}

We now introduce the notation required for Theorem \ref{thm:Bounds on the Resolvent Degree of the Sporadic Groups}. 

\begin{notation}[Notation for Theorem \ref{thm:Bounds on the Resolvent Degree of the Sporadic Groups}]\label{not:Notation for Bounds on the Resolvent Degree of the Sporadic Groups}
Let $G$ be a sporadic group. When $m_d\lp \rho_G \rp = j$, we denote a basis for $\Sym_{a.G}^d\lp V_G^\vee \rp$ by $f_{d,1}^G,\dotsc,f_{d,j}^G$. When $j=1$, we simply write $f_d^G$. Without loss of generality, we order the basis such that the basis elements which are algebraically independent from lower degree invariants are listed first; see Remark \ref{rem:Non-Uniqueness of XG, YG, ZG} for more details. When the group $G$ is clear from context, we omit the superscript. Using this notation, we now define the relevant $G$-invariant subvarieties of $V_G$. Specifically, we will define a variety $X_G = Z_G \cap Y_G$. We begin with the cases where $Y_G$ is non-trivial:

\begin{align*}
    Z_{\MOO} &= \VV\lp f_{4,1} \rp, &Y_{\MOO} &= \VV\lp f_2,f_3 \rp, &X_{\MOO} &= \VV\lp f_2,f_3,f_{4,1} \rp \subseteq \PP\lp V_{\MOO} \rp = \PP^9,\\
	Z_{\MOW} &= \VV\lp f_{4,1} \rp, &Y_{\MOW} &= \VV\lp f_2,f_3 \rp, &X_{\MOW} &= \VV\lp f_2,f_3,f_{4,1} \rp \subseteq \PP\lp V_{\MOW} \rp = \PP^{10},\\
    Z_{\MWH} &= \VV\lp f_{4,1}, f_{5,1} \rp, &Y_{\MWH} &= \VV\lp f_2, f_3 \rp, &X_{\MWH} &= \VV\lp f_2, f_3, f_{4,1}, f_{5,1} \rp \subseteq \PP\lp V_{\MWH} \rp = \PP^{21}, \\
    Z_{\HS} &= \VV\lp f_{4,1}, f_5 \rp, &Y_{\HS} &= \VV\lp f_2 \rp, &X_{\HS} &= \VV\lp f_2, f_{4,1}, f_5 \rp \subseteq \PP\lp V_{\HS} \rp = \PP^{21},\\
	Z_{\MCL} &= \VV\lp f_5 \rp, &Y_{\MCL} &= \VV\lp f_2 \rp, &X_{\MCL} &= \VV\lp f_2, f_5 \rp \subseteq \PP\lp V_{\MCL} \rp = \PP^{21},\\
	Z_{\MWF} &= \VV\lp f_{4,1}, f_{5,1} \rp, &Y_{\MWF} &= \VV\lp f_2, f_3 \rp, &X_{\MWF} &= \VV\lp f_2, f_3, f_{4,1}, f_{5,1} \rp \subseteq \PP\lp V_{\MWF} \rp = \PP^{22}, \\
    Z_{\COH} &= \VV\lp f_{6,1} \rp, &Y_{\COH} &= \VV\lp f_2 \rp, &X_{\COH} &= \VV\lp f_2,f_{6,1} \rp \subseteq \PP\lp V_{\COH} \rp = \PP^{22},\\
	Z_{\COW} &= \VV\lp f_{8,1} \rp, &Y_{\COW} &= \VV\lp f_2 \rp, &X_{\COW} &= \VV\lp f_2, f_{8,1} \rp \subseteq \PP\lp V_{\COW} \rp = \PP^{22},\\
	Z_{\COO} &= \VV\lp f_{12,1} \rp, &Y_{\COO} &= \VV\lp f_2 \rp, &X_{\COO} &= \VV\lp f_2, f_{12,1} \rp \subseteq \PP\lp V_{\COO} \rp = \PP^{23},\\
    Z_{\JO} &= \VV\lp f_3, f_{4,1}, f_{4,2} \rp, &Y_{\JO} &= \VV\lp f_2 \rp, &X_{\JO} &= \VV\lp f_2, f_3, f_{4,1}, f_{4,2} \rp \subseteq \PP\lp V_{\JO} \rp = \PP^{55},
\end{align*}

\begin{align*}
    Z_{\FIWW} &= \VV\lp f_6, f_{8,1} \rp, &Y_{\FIWW} &= \VV\lp f_2 \rp, &X_{\FIWW} &= \VV\lp f_2, f_6, f_{8,1} \rp \subseteq \PP\lp V_{\FIWW} \rp = \PP^{77},\\
	Z_{\HN} &= \VV\lp f_6, f_7 \rp, &Y_{\HN} &= \VV\lp f_2 \rp, &X_{\HN} &= \VV\lp f_2, f_6, f_7 \rp \subseteq \PP\lp V_{\HN} \rp = \PP^{132},\\
	Z_{\THO} &= \VV\lp f_{8,1}, f_{8,2} \rp, &Y_{\THO} &= \VV\lp f_2 \rp, &X_{\THO} &= \VV\lp f_2, f_{8,1}, f_{8,2} \rp \subseteq \PP\lp V_{\THO} \rp = \PP^{247}.
\end{align*}

\noindent
In the following cases, we have $X_G = Z_G$ and $Y_G = \PP\lp V_G \rp$:
\begin{align*}
    Z_{\MWW} = X_{\MWW} &= \VV\lp f_4 \rp  \subseteq \PP\lp V_{\MWW} \rp = \PP^{9}, \\
	Z_{\SUZ} = X_{\SUZ} &= \VV\lp f_{12} \rp \subseteq \PP\lp V_{\SUZ} \rp = \PP^{11},\\
	Z_{\JH} = X_{\JH} &= \VV\lp f_6 \rp \subseteq \PP\lp V_{\JH} \rp = \PP^{17},\\
    Z_{\RU} = X_{\RU} &= \VV\lp f_4 \rp \subseteq \PP\lp V_{\RU} \rp = \PP^{27},\\
	Z_{\HE} = X_{\HE} &= \VV\lp f_3, f_4 \rp \subseteq \PP\lp V_{\HE} \rp = \PP^{50},\\
	Z_{\ON} = X_{\ON} &= \VV\lp f_{6,1}, f_{6,2}, f_{6,3} \rp \subseteq \PP\lp V_{\ON} \rp = \PP^{341},\\
    Z_{\FIWH} = X_{\FIWH} &= \VV\lp f_2, f_3, f_{4,1}, f_{5,1}, f_{5,2} \rp \subseteq \PP\lp V_{\FIWH} \rp = \PP^{781},\\
	Z_{\FIWF} = X_{\FIWF} &= \VV\lp f_3, f_{6,1}, f_{6,2} \rp \subseteq \PP\lp V_{\FIWF} \rp = \PP^{782}, \\
    Z_{\JF} = X_{\JF} &= \VV\lp f_4, f_{6,1}, f_{6,2}, f_{7,1} \rp \subseteq \PP\lp V_{\JF} \rp = \PP^{1332}, \\
	Z_{\LY} = X_{\LY} &= \VV\lp f_{6,1}, f_{6,2}, f_{6,3}, f_{6,4} \rp \subseteq \PP\lp V_{\LY} \rp = \PP^{2479},\\
	Z_{\B} = X_{\B} &= \VV\lp f_2, f_{4,1}, f_{6,1}, f_{8,1}, f_{8,2} \rp \subseteq \PP\lp V_{\B} \rp = \PP^{4370},\\
	Z_{\MO} = X_{\MO} &= \VV\lp f_2, f_3, f_{4,1}, f_{5,1}, f_{6,1}, f_{6,2}, f_{6,3}, f_{7,1} \rp \subseteq \PP\lp V_{\MO} \rp = \PP^{196882}.
\end{align*}

\noindent
In the exceptional case of $\JW$, we further have that $X_{\JW} = Z_{\JW} = Y_{\JW} = \PP\lp V_{\JW} \rp = \PP^5$. 
\end{notation}

\begin{remark}[Non-Uniqueness of $X_G$, $Y_G$, $Z_G$]\label{rem:Non-Uniqueness of XG, YG, ZG}
We note that the $X_G$, $Y_G,$ and $Z_G$ are only defined up to a choice of invariant polynomials. As an example, consider the case $G=\COH$. First, note that $m_2\lp \rho_{\COH} \rp=1$ and thus $f_2$ is unique up to a constant. Next, $m_4\lp \rho_{\COH} \rp=1$ and thus the only degree four invariant, up to scaling, is $\lp f_2 \rp^2$. Finally, $m_6\lp \rho_{\COH} \rp = 2$ and so we take $f_{6,1}$ to be any polynomial in $\Sym_{\COH}^6\lp V_{\COH}^\vee \rp \setminus \SPAN\lb \lp f_2 \rp^3 \rb$ and $f_{6,2}$ to be any polynomial in $\SPAN\lb \lp f_2 \rp^3 \rb$. However, all of the arguments that follow depend only on the degrees of the intersections defining $X_G$, $Y_G,$ and $Z_G$ and are consequently independent of this choice. 
\end{remark}

We are now ready to give upper bounds on the resolvent degree of the sporadic groups.

\begin{theorem}[Bounds on Resolvent Degree of the Sporadic Groups]\label{thm:Bounds on the Resolvent Degree of the Sporadic Groups}
For each sporadic group $G$, we have
\begin{equation*}
	\RD_k(G) \leq \dim_{\CC}(X_G).
\end{equation*}
\noindent
for every field $k$. Further, for $G$ not equal to $\MOO,\MOW,\MWH,\MWF$, the variety $X_G$ is $\RD_{\CC}^{\leq d_G}$-versal for $d_G = \RDC(\deg(Z_G))$.
\end{theorem}

\begin{remark}
    We expect that the variety $X_G$ is $\RD_k^{\leq d_G}$-versal for the Mathieu groups $\MOO,\MOW,\MWH,\MWF$, but proving this requires new techniques.
\end{remark}

\noindent
We will also give an explicit form of Theorem \ref{thm:Bounds on the Resolvent Degree of the Sporadic Groups} in Corollary \ref{cor:Explicit Form of Bounds on the Resolvent Degree of the Sporadic Groups}. 

\begin{proof}
As observed above, by \cite[Theorems 1.2 and 1.3]{Reichstein2022}, it suffices to prove the theorem for $k=\CC$. When $G$ is one of $\MOO$, $\MOW$, $\MWH$, or $\MWF$, the upper bounds on $\RD_{\CC}(G)$ are classical. In these cases, $M_n < S_n$ and the bounds follow from the inequality $\RD_{\CC}(M_n)\le \RD_{\CC}(S_n)$ \cite[Lemma 3.13]{FarbWolfson2019} and the classical upper bounds on $\RD_{\CC}(S_n)$ (see \cite[Theorem 3.7]{Sutherland2021} for the construction and the bounds on $S_{23}, S_{24}$; see \cite{FarbWolfson2019,Sutherland2021} for modern references for $S_{11}, S_{12}$ or \cite{Hilbert1927,Segre1945} for classical references). 

We now restrict to the case where $G$ is not one of $\MOO$, $\MOW$, $\MWH$, $\MWF$. Note that $d_G \leq \dim_{\CC}\lp X_G \rp$, hence Proposition \ref{prop:RDK(G) via RDKD-versality} yields that we need only show that each $X_G$ is $\RD_{\CC}^{\leq d_G}$-versal. Since finite groups are smooth, Theorem \ref{thm:RDKD-versality and Special Points} allows us to reduce to showing that a) $X_G$ is generically free, and b) for every $G$-torsor $T \ra \spec(K)$ with $K$ finitely generated over $\CC$, $K^{\lp d_G\rp}$-points are dense in $\TX_G$.

We start with generic freeness.  Since $G$ is simple, and the representation $\rho_{G}$ is irreducible, we see that $\PP(V_G)$ has no fixed points and thus is a faithful $G$-variety. By Lemma~\ref{lem:faithful and irred suffices}, it suffices to show that $X_G$ is irreducible, but this follows for degree reasons.  Indeed, for all simple sporadic $G$ not equal to $\MOO$, $\MOW$, $\MWH$, or $\MWF$, the degree of $X_G$ is less than the cardinality of the smallest permutation representation of $G$ (see Figure~\ref{table:Dimensions of Projective Representations and Degrees of Invariants} and preceeding discussion).  Since $X_G$ has at most $\deg X_G$ irreducible components, and since $G$ permutes them, we conclude that $X_G$ is irreducible, and thus generically free.  

To apply Theorem~\ref{thm:RDKD-versality and Special Points} to conclude the $\RD_{\CC}^{\leq d_G}$-versality of $X_G$, we just need to show that $K^{\lp d_G\rp}$ points are dense in every twisted form of $X_G$.  As observed in \cite[Proof of Proposition 14.1 (p.33)]{Reichstein2022}, the $G$-equivariant closed immersion $X_G \injects \PP\lp V_G \rp$ naturally induces the closed immersion $\TXG \injects \prescript{T}{}{\PP\lp V_G \rp}$. Note that $\prescript{T}{}{\PP\lp V_G \rp}$ is a Severi--Brauer variety and thus splits over $\KS \subseteq K^{\lp d_G \rp}$ by the Merkurjev-Suslin theorem \cite{MerkurjevSuslin1983}. It follows that $\TXG$ is an intersection of hypersurfaces in $\PP\lp V_G \rp$ over $K^{\lp d_G \rp}$ of the same degrees as for $X_G$. Indeed, the same argument applies to $\TYG \injects \prescript{T}{}{\PP\lp V_G \rp}$ and $\TZG \injects \prescript{T}{}{\PP\lp V_G \rp}$ and we have $\TXG = \TZG \cap \TYG$ over $K^{\lp d_G \rp}$.

Now, observe that for each $G$, we either have that $Y_G = \PP\lp V_G \rp$ (and thus $X_G = Z_G$) or $Y_G = \VV\lp f_2^G \rp$. In the first case, we have $\deg\lp X_G \rp = \deg\lp Z_G \rp = d_G$, and the density of $K^{\lp d_G \rp}$-points is immediate. In the second case, $Y_G$ is a quadric hypersurface and $X_G = Y_G \cap Z_G$, hence \cite[Lemma 14.4(b)]{Reichstein2022} yields that $K^{\lp d_G \rp}$-points are dense on $X_G$.
\end{proof}

\begin{corollary}[Explicit Form of Theorem \ref{thm:Bounds on the Resolvent Degree of the Sporadic Groups}]\label{cor:Explicit Form of Bounds on the Resolvent Degree of the Sporadic Groups}
For any field $k$, we have
\begin{align*}
	\RD_k(\JW) &\leq 5, &\RD_k(\MWF) &\leq 18, &\RD_k(\HE) &\leq 48, &\RD_k(\FIWH) &\leq 776,\\
	\RD_k(\MOO) &\leq 6, &\RD_k(\HS) &\leq 18, &\RD_k(\JO) &\leq 51, &\RD_k(\FIWF) &\leq 779,\\
	\RD_k(\MOW) &\leq 7, &\RD_k(\MCL) &\leq 19, &\RD_k(\FIWW) &\leq 74, &\RD_k(\JF) &\leq 1328,\\
	\RD_k(\MWW) &\leq 8, &\RD_k(\COH) &\leq 20, &\RD_k(\HN) &\leq 129, &\RD_k(\LY) &\leq 2475,\\
	\RD_k(\SUZ) &\leq 10, &\RD_k(\COW) &\leq 20, &\RD_k(\THO) &\leq 244, &\RD_k(\B) &\leq 4365,\\
	\RD_k(\JH) &\leq 16, &\RD_k(\COO) &\leq 21, &\RD_k(\ON) &\leq 338, &\RD_k(\MO) &\leq 196874.\\
	\RD_k(\MWH) &\leq 17, &\RD_k(\RU) &\leq 26, 
\end{align*}
\end{corollary}

\begin{remark}[Further Expectations]\label{rem:Further Expectations}
In the cases where $G$ is one of $\MWW$, $\RU$, $\HE$, $\FIWH$, or $\FIWF$, we expect that we can do slightly better. Indeed, we believe that we can replace $X_G$, $Y_G$, and $Z_G$ with
\begin{align*}
    \wt{Z}_{\MWW} &= \VV\lp f_6 \rp, &\wt{Y}_{\MWW} &= \VV\lp f_4 \rp, &\wt{X}_{\MWW} &= \VV\lp f_4, f_6 \rp \subseteq \PP\lp V_{\MWW} \rp = \PP^{9},\\	
    \wt{Z}_{\RU} &= \VV\lp f_8 \rp, &\wt{Y}_{\RU} &= \VV\lp f_4 \rp, &\wt{X}_{\RU} &= \VV\lp f_4, f_8 \rp \subseteq \PP\lp V_{\RU} \rp = \PP^{27},\\
	\wt{Z}_{\HE} &= \VV\lp f_4, f_5 \rp, &\wt{Y}_{\HE} &= \VV\lp f_3 \rp, &\wt{X}_{\HE} &= \VV\lp f_3, f_4, f_5 \rp \subseteq \PP\lp V_{\HE} \rp = \PP^{50},\\
	\wt{Z}_{\FIWH} &= \VV\lp f_4, f_{5,1}, f_{5,2}, f_6 \rp, &\wt{Y}_{\FIWH} &= \VV\lp f_2, f_3 \rp, &\wt{X}_{\FIWH} &= \VV\lp f_2, f_3, f_4, f_{5,1}, f_{5,2}, f_6 \rp \subseteq \PP\lp V_{\FIWH} \rp = \PP^{781},\\
	\wt{Z}_{\FIWF} &= \VV\lp f_{6,1}, f_{6,2}, f_9 \rp, &\wt{Y}_{\FIWF} &= \VV\lp f_3 \rp, &\wt{X}_{\FIWF} &= \VV\lp f_3, f_{6,1}, f_{6,2}, f_9 \rp \subseteq \PP\lp V_{\FIWF} \rp = \PP^{782}.
\end{align*}

\noindent
In these cases, one can use the polar cone methods of \cite{Sutherland2021} to construct linear subvarieties of suitable dimension and satisfactorily low resolvent degree on each $\wt{Y}_G$. However, new methods are required to show that each $\wt{X}_G = \wt{Y}_G \cap \wt{Z}_G$ is generically free. 
\end{remark}

\begin{remark}[The Unique Case of $\MOW$, II]\label{rem:The Unique Case of MOW, II}
    As noted in Remark \ref{rem:The Unique Case of MOW}, $\rho_{\MOW}$ is the only case where we are not using a projective representation of minimal dimension. As we have seen, we constructed
    \begin{equation*}
        X_{\MOW} = \VV\lp f_2, f_3, f_4 \rp \subseteq \PP^{10} = \PP\lp V_{\MOW} \rp.
    \end{equation*}

    Now, let $\wt{\rho}_{\MOW}$ be the first 10-dimensional representation of the Schur cover $2.\MOW$. The lowest degree invariants of $\wt{\rho}_{\MOW}$ have degrees 6, 8, 8, and 8 respectively. Consequently, the best analogous construction would yield
    \begin{equation*}
        \wt{X}_{\MOW} = \VV\lp f_6 \rp \subseteq \PP^9 = \PP\lp \wt{V}_{\MOW} \rp.
    \end{equation*}

    \noindent
    Since $\dim\lp X_{\MOW} \rp < \dim\lp \wt{X}_{\MOW} \rp$, $X_{\MOW}$ is the preferred construction. 
    
    For every other sporadic group $G$ with non-trivial Schur cover $a.G$, either there is not a lower dimensional projective representation or the dimension of the new projective representation is small enough to outweigh any differences in the invariant theory.
\end{remark}

\subsection{Context for Sporadic Group Bounds}\label{subsec:Context for Sporadic Group Bounds}
We conclude by providing additional context for these numerical results and connect the bounds for the Mathieu groups to known bounds for symmetric groups.

\paragraph{Mathieu Groups and Symmetric Groups}\label{par:Mathieu Groups and Symmetric Groups}

As noted in Section \ref{sec:Introduction}, $\RD_k(G)$ for finite simple groups has only been addressed in the literature when $G$ is a cyclic group (for which $\RD_k \equiv 1$ by Kummer theory), 
an alternating group \cite{HeberleSutherland2023,Sutherland2021,Wolfson2021}, when $G = W\lp E_6\rp^+$, $W\lp E_7 \rp^+$, or $W\lp E_8 \rp^+$ (see \cite[Section 8]{FarbWolfson2019} for $E_6$ and $E_7$, see \cite[Proposition 15.1]{Reichstein2022} for $E_6$, $E_7$, and $E_8$), or when $G=\PSL(2,7)$ (\cite{Klein1879}, \cite[Proposition 4.13]{FarbKisinWolfson2023}). Nonetheless, each of the Mathieu groups have explicit embeddings $M_n \injects S_n$ for $n=11$, 12, 22, 23, 24 and thus $\RD_k(M_n) \leq \RD_k(S_n)$ \cite[Lemma 3.13]{FarbWolfson2019}. At present, this paper contributes nothing new for $\MOO,\MOW,\MWH,\MWF$. In the case of $\MWW$, our bound of 8 significantly beats the bound $\RD_\CC(S_{22}) \leq 16$ \cite[Theorem 3.7]{Sutherland2021}. As remarked above, it would be interesting to confirm that for $\MOO,\MOW,\MWH,\MWF$, the variety $X_G$ is also $\RD_{\CC}^{\leq d_G}$-versal.  It would also be interesting to know for which $n$ and which fields $k$ we have $\RD_k(M_n) < \RD_k(S_n)$, and for which $n$ and which fields $k$ we have $\RD_k(M_n) = \RD_k(S_n)$.

\paragraph{Relations Between the Sporadic Groups}\label{par::Relations Between the Sporadic Groups}

For a finite group $H$ and a subgroup $H'$, we have $\RD_k(H') \leq \RD_k(H)$ \cite[Lemma 3.13]{FarbWolfson2019}. Additionally, for a short exact sequence of algebraic groups
\begin{equation*}
	1 \ra A \ra B \ra C \ra 1,
\end{equation*}

\noindent
we have that $\RD_k(B) \leq \max \lb \RD_k(A), \RD_k(C) \rb$ (\cite[Theorem 3.3]{FarbWolfson2019} for finite groups, \cite[Proposition 10.8]{Reichstein2022} in general).

For any sporadic group $S$ which is a subquotient of another sporadic group $G$, with $H < G$ and $S = H/H'$, we only have the inequalities
\begin{align*}
 \RD_k(H) &\leq \RD_k(G), &\RD_k(H) &\leq \max\lb \RD_k(S), \RD_k(H') \rb.
\end{align*}

\noindent
Nonetheless, it is natural to ask how the bounds given by $\dim\lp X_S \rp$ and $\dim\lp X_G \rp$ in Theorem \ref{thm:Bounds on the Resolvent Degree of the Sporadic Groups} compare.

The $\ATLAS$ contains a table of which sporadic groups $S$ are subquotients of another sporadic group $G$ \cite[p.238]{ConwayCurtisNortonParkerWilson1985}, which we include below. Let $(G,S)$ be the cell corresponding to row $G$ and column $S$. Note that $(G,S)$ has a $+$ with a green background if $S$ is a subquotient of $G$; $(G,S)$ has a $-$ with a red background if $S$ is not a subquotient of $G$; $(G,S)$ has a $\bullet$ with a yellow background when $G=S$; and $(G,S)$ is black when $|G| < |S|$.

We implement two changes from the table in the $\ATLAS$. Firstly, at the time of publishing the $\ATLAS$, it was unknown if $\JO$ is a subquotient of $\MO$. Wilson showed that $\JO$ is not a subgroup of $\MO$ in \cite{Wilson1986}, which completed the proof that $\JO$ is not a subquotient of $\MO$, and we have updated the $\lp \MO, \JO \rp$ cell as a result. Secondly, due to page size restrictions, we split the single table in the $\ATLAS$ into two smaller tables: Figures \ref{table:Sporadic Groups Table from the ATLAS} and \ref{table:Sporadic Groups Table from the ATLAS, II}.

Finally, we note that
\begin{align*}
	&\dim\lp X_{\JW} \rp < \dim\lp X_{\MOO} \rp \leq \dim\lp X_{\MOW} \rp < \dim\lp X_{\MWW} \rp < \dim\lp X_{\SUZ} \rp < \dim\lp X_{\JH} \rp < \dim\lp X_{\MWH} \rp\\
	< \ &\dim\lp X_{\MWF} \rp = \dim\lp X_{\HS} \rp < \dim\lp X_{\MCL} \rp < \dim\lp X_{\COH} \rp = \dim\lp X_{\COW} \rp < \dim\lp X_{\COO} \rp < \dim\lp X_{\RU} \rp\\
	< \ &\dim\lp X_{\HE} \rp < \dim\lp X_{\JO} \rp < \dim\lp X_{\FIWW} \rp < \dim\lp X_{\HN} \rp < \dim\lp X_{\THO} \rp < \dim\lp X_{\ON} \rp < \dim\lp X_{\FIWH} \rp\\
	< \ &\dim\lp X_{\FIWF} \rp < \dim\lp X_{\JF} \rp < \dim\lp X_{\LY} \rp < \dim\lp X_{\B} \rp < \dim\lp X_{\MO} \rp,
\end{align*}

\noindent
and thus one can verify using Figures \ref{table:Sporadic Groups Table from the ATLAS} and \ref{table:Sporadic Groups Table from the ATLAS, II} that whenever $S$ is a subquotient of $G$, we have $\dim\lp X_S \rp \leq \dim\lp X_G \rp$.

\begin{figure}[ht!]
\begin{center}
\caption{Sporadic Groups Subquotient Table from the $\ATLAS$, I}
\label{table:Sporadic Groups Table from the ATLAS}
\begin{tabular}{ccccccccccccccccc}
        & $\MOO$                            & $\MOW$                            & $\JO$                             & $\MWW$                            & $\JW$                             & $\MWH$                            & $\HS$                             & $\JH$                             & $\MWF$                            & $\MCL$                            & $\HE$                             & $\RU$                             & $\SUZ$                            & $\ON$                             & $\COH$                            & $\COW$                            \\
$\MOO$  & \cellcolor[HTML]{FFFC9E}$\bullet$ & \cellcolor[HTML]{000000}          & \cellcolor[HTML]{000000}          & \cellcolor[HTML]{000000}          & \cellcolor[HTML]{000000}          & \cellcolor[HTML]{000000}          & \cellcolor[HTML]{000000}          & \cellcolor[HTML]{000000}          & \cellcolor[HTML]{000000}          & \cellcolor[HTML]{000000}          & \cellcolor[HTML]{000000}          & \cellcolor[HTML]{000000}          & \cellcolor[HTML]{000000}          & \cellcolor[HTML]{000000}          & \cellcolor[HTML]{000000}          & \cellcolor[HTML]{000000}          \\
$\MOW$  & \cellcolor[HTML]{34FF34}$+$       & \cellcolor[HTML]{FFFC9E}$\bullet$ & \cellcolor[HTML]{000000}          & \cellcolor[HTML]{000000}          & \cellcolor[HTML]{000000}          & \cellcolor[HTML]{000000}          & \cellcolor[HTML]{000000}          & \cellcolor[HTML]{000000}          & \cellcolor[HTML]{000000}          & \cellcolor[HTML]{000000}          & \cellcolor[HTML]{000000}          & \cellcolor[HTML]{000000}          & \cellcolor[HTML]{000000}          & \cellcolor[HTML]{000000}          & \cellcolor[HTML]{000000}          & \cellcolor[HTML]{000000}          \\
$\JO$   & \cellcolor[HTML]{FD6864}$-$       & \cellcolor[HTML]{FD6864}$-$       & \cellcolor[HTML]{FFFC9E}$\bullet$ & \cellcolor[HTML]{000000}          & \cellcolor[HTML]{000000}          & \cellcolor[HTML]{000000}          & \cellcolor[HTML]{000000}          & \cellcolor[HTML]{000000}          & \cellcolor[HTML]{000000}          & \cellcolor[HTML]{000000}          & \cellcolor[HTML]{000000}          & \cellcolor[HTML]{000000}          & \cellcolor[HTML]{000000}          & \cellcolor[HTML]{000000}          & \cellcolor[HTML]{000000}          & \cellcolor[HTML]{000000}          \\
$\MWW$  & \cellcolor[HTML]{FD6864}$-$       & \cellcolor[HTML]{FD6864}$-$       & \cellcolor[HTML]{FD6864}$-$       & \cellcolor[HTML]{FFFC9E}$\bullet$ & \cellcolor[HTML]{000000}          & \cellcolor[HTML]{000000}          & \cellcolor[HTML]{000000}          & \cellcolor[HTML]{000000}          & \cellcolor[HTML]{000000}          & \cellcolor[HTML]{000000}          & \cellcolor[HTML]{000000}          & \cellcolor[HTML]{000000}          & \cellcolor[HTML]{000000}          & \cellcolor[HTML]{000000}          & \cellcolor[HTML]{000000}          & \cellcolor[HTML]{000000}          \\
$\JW$   & \cellcolor[HTML]{FD6864}$-$       & \cellcolor[HTML]{FD6864}$-$       & \cellcolor[HTML]{FD6864}$-$       & \cellcolor[HTML]{FD6864}$-$       & \cellcolor[HTML]{FFFC9E}$\bullet$ & \cellcolor[HTML]{000000}          & \cellcolor[HTML]{000000}          & \cellcolor[HTML]{000000}          & \cellcolor[HTML]{000000}          & \cellcolor[HTML]{000000}          & \cellcolor[HTML]{000000}          & \cellcolor[HTML]{000000}          & \cellcolor[HTML]{000000}          & \cellcolor[HTML]{000000}          & \cellcolor[HTML]{000000}          & \cellcolor[HTML]{000000}          \\
$\MWH$  & \cellcolor[HTML]{34FF34}$+$       & \cellcolor[HTML]{FD6864}$-$       & \cellcolor[HTML]{FD6864}$-$       & \cellcolor[HTML]{34FF34}$+$       & \cellcolor[HTML]{FD6864}$-$       & \cellcolor[HTML]{FFFC9E}$\bullet$ & \cellcolor[HTML]{000000}          & \cellcolor[HTML]{000000}          & \cellcolor[HTML]{000000}          & \cellcolor[HTML]{000000}          & \cellcolor[HTML]{000000}          & \cellcolor[HTML]{000000}          & \cellcolor[HTML]{000000}          & \cellcolor[HTML]{000000}          & \cellcolor[HTML]{000000}          & \cellcolor[HTML]{000000}          \\
$\HS$   & \cellcolor[HTML]{34FF34}$+$       & \cellcolor[HTML]{FD6864}$-$       & \cellcolor[HTML]{FD6864}$-$       & \cellcolor[HTML]{34FF34}$+$       & \cellcolor[HTML]{FD6864}$-$       & \cellcolor[HTML]{FD6864}$-$       & \cellcolor[HTML]{FFFC9E}$\bullet$ & \cellcolor[HTML]{000000}          & \cellcolor[HTML]{000000}          & \cellcolor[HTML]{000000}          & \cellcolor[HTML]{000000}          & \cellcolor[HTML]{000000}          & \cellcolor[HTML]{000000}          & \cellcolor[HTML]{000000}          & \cellcolor[HTML]{000000}          & \cellcolor[HTML]{000000}          \\
$\JH$   & \cellcolor[HTML]{FD6864}$-$       & \cellcolor[HTML]{FD6864}$-$       & \cellcolor[HTML]{FD6864}$-$       & \cellcolor[HTML]{FD6864}$-$       & \cellcolor[HTML]{FD6864}$-$       & \cellcolor[HTML]{FD6864}$-$       & \cellcolor[HTML]{FD6864}$-$       & \cellcolor[HTML]{FFFC9E}$\bullet$ & \cellcolor[HTML]{000000}          & \cellcolor[HTML]{000000}          & \cellcolor[HTML]{000000}          & \cellcolor[HTML]{000000}          & \cellcolor[HTML]{000000}          & \cellcolor[HTML]{000000}          & \cellcolor[HTML]{000000}          & \cellcolor[HTML]{000000}          \\
$\MWF$  & \cellcolor[HTML]{34FF34}$+$       & \cellcolor[HTML]{34FF34}$+$       & \cellcolor[HTML]{FD6864}$-$       & \cellcolor[HTML]{34FF34}$+$       & \cellcolor[HTML]{FD6864}$-$       & \cellcolor[HTML]{34FF34}$+$       & \cellcolor[HTML]{FD6864}$-$       & \cellcolor[HTML]{FD6864}$-$       & \cellcolor[HTML]{FFFC9E}$\bullet$ & \cellcolor[HTML]{000000}          & \cellcolor[HTML]{000000}          & \cellcolor[HTML]{000000}          & \cellcolor[HTML]{000000}          & \cellcolor[HTML]{000000}          & \cellcolor[HTML]{000000}          & \cellcolor[HTML]{000000}          \\
$\MCL$  & \cellcolor[HTML]{34FF34}$+$       & \cellcolor[HTML]{FD6864}$-$       & \cellcolor[HTML]{FD6864}$-$       & \cellcolor[HTML]{34FF34}$+$       & \cellcolor[HTML]{FD6864}$-$       & \cellcolor[HTML]{FD6864}$-$       & \cellcolor[HTML]{FD6864}$-$       & \cellcolor[HTML]{FD6864}$-$       & \cellcolor[HTML]{FD6864}$-$       & \cellcolor[HTML]{FFFC9E}$\bullet$ & \cellcolor[HTML]{000000}          & \cellcolor[HTML]{000000}          & \cellcolor[HTML]{000000}          & \cellcolor[HTML]{000000}          & \cellcolor[HTML]{000000}          & \cellcolor[HTML]{000000}          \\
$\HE$   & \cellcolor[HTML]{FD6864}$-$       & \cellcolor[HTML]{FD6864}$-$       & \cellcolor[HTML]{FD6864}$-$       & \cellcolor[HTML]{FD6864}$-$       & \cellcolor[HTML]{FD6864}$-$       & \cellcolor[HTML]{FD6864}$-$       & \cellcolor[HTML]{FD6864}$-$       & \cellcolor[HTML]{FD6864}$-$       & \cellcolor[HTML]{FD6864}$-$       & \cellcolor[HTML]{FD6864}$-$       & \cellcolor[HTML]{FFFC9E}$\bullet$ & \cellcolor[HTML]{000000}          & \cellcolor[HTML]{000000}          & \cellcolor[HTML]{000000}          & \cellcolor[HTML]{000000}          & \cellcolor[HTML]{000000}          \\
$\RU$   & \cellcolor[HTML]{FD6864}$-$       & \cellcolor[HTML]{FD6864}$-$       & \cellcolor[HTML]{FD6864}$-$       & \cellcolor[HTML]{FD6864}$-$       & \cellcolor[HTML]{FD6864}$-$       & \cellcolor[HTML]{FD6864}$-$       & \cellcolor[HTML]{FD6864}$-$       & \cellcolor[HTML]{FD6864}$-$       & \cellcolor[HTML]{FD6864}$-$       & \cellcolor[HTML]{FD6864}$-$       & \cellcolor[HTML]{FD6864}$-$       & \cellcolor[HTML]{FFFC9E}$\bullet$ & \cellcolor[HTML]{000000}          & \cellcolor[HTML]{000000}          & \cellcolor[HTML]{000000}          & \cellcolor[HTML]{000000}          \\
$\SUZ$  & \cellcolor[HTML]{34FF34}$+$       & \cellcolor[HTML]{34FF34}$+$       & \cellcolor[HTML]{FD6864}$-$       & \cellcolor[HTML]{FD6864}$-$       & \cellcolor[HTML]{34FF34}$+$       & \cellcolor[HTML]{FD6864}$-$       & \cellcolor[HTML]{FD6864}$-$       & \cellcolor[HTML]{FD6864}$-$       & \cellcolor[HTML]{FD6864}$-$       & \cellcolor[HTML]{FD6864}$-$       & \cellcolor[HTML]{FD6864}$-$       & \cellcolor[HTML]{FD6864}$-$       & \cellcolor[HTML]{FFFC9E}$\bullet$ & \cellcolor[HTML]{000000}          & \cellcolor[HTML]{000000}          & \cellcolor[HTML]{000000}          \\
$\ON$   & \cellcolor[HTML]{34FF34}$+$       & \cellcolor[HTML]{FD6864}$-$       & \cellcolor[HTML]{34FF34}$+$       & \cellcolor[HTML]{FD6864}$-$       & \cellcolor[HTML]{FD6864}$-$       & \cellcolor[HTML]{FD6864}$-$       & \cellcolor[HTML]{FD6864}$-$       & \cellcolor[HTML]{FD6864}$-$       & \cellcolor[HTML]{FD6864}$-$       & \cellcolor[HTML]{FD6864}$-$       & \cellcolor[HTML]{FD6864}$-$       & \cellcolor[HTML]{FD6864}$-$       & \cellcolor[HTML]{FD6864}$-$       & \cellcolor[HTML]{FFFC9E}$\bullet$ & \cellcolor[HTML]{000000}          & \cellcolor[HTML]{000000}          \\
$\COH$  & \cellcolor[HTML]{34FF34}$+$       & \cellcolor[HTML]{34FF34}$+$       & \cellcolor[HTML]{FD6864}$-$       & \cellcolor[HTML]{34FF34}$+$       & \cellcolor[HTML]{FD6864}$-$       & \cellcolor[HTML]{34FF34}$+$       & \cellcolor[HTML]{34FF34}$+$       & \cellcolor[HTML]{FD6864}$-$       & \cellcolor[HTML]{FD6864}$-$       & \cellcolor[HTML]{34FF34}$+$       & \cellcolor[HTML]{FD6864}$-$       & \cellcolor[HTML]{FD6864}$-$       & \cellcolor[HTML]{FD6864}$-$       & \cellcolor[HTML]{FD6864}$-$       & \cellcolor[HTML]{FFFC9E}$\bullet$ & \cellcolor[HTML]{000000}          \\
$\COW$  & \cellcolor[HTML]{34FF34}$+$       & \cellcolor[HTML]{FD6864}$-$       & \cellcolor[HTML]{FD6864}$-$       & \cellcolor[HTML]{34FF34}$+$       & \cellcolor[HTML]{FD6864}$-$       & \cellcolor[HTML]{34FF34}$+$       & \cellcolor[HTML]{34FF34}$+$       & \cellcolor[HTML]{FD6864}$-$       & \cellcolor[HTML]{FD6864}$-$       & \cellcolor[HTML]{34FF34}$+$       & \cellcolor[HTML]{FD6864}$-$       & \cellcolor[HTML]{FD6864}$-$       & \cellcolor[HTML]{FD6864}$-$       & \cellcolor[HTML]{FD6864}$-$       & \cellcolor[HTML]{FD6864}$-$       & \cellcolor[HTML]{FFFC9E}$\bullet$ \\
$\FIWW$ & \cellcolor[HTML]{34FF34}$+$       & \cellcolor[HTML]{34FF34}$+$       & \cellcolor[HTML]{FD6864}$-$       & \cellcolor[HTML]{34FF34}$+$       & \cellcolor[HTML]{FD6864}$-$       & \cellcolor[HTML]{FD6864}$-$       & \cellcolor[HTML]{FD6864}$-$       & \cellcolor[HTML]{FD6864}$-$       & \cellcolor[HTML]{FD6864}$-$       & \cellcolor[HTML]{FD6864}$-$       & \cellcolor[HTML]{FD6864}$-$       & \cellcolor[HTML]{FD6864}$-$       & \cellcolor[HTML]{FD6864}$-$       & \cellcolor[HTML]{FD6864}$-$       & \cellcolor[HTML]{FD6864}$-$       & \cellcolor[HTML]{FD6864}$-$       \\
$\HN$   & \cellcolor[HTML]{34FF34}$+$       & \cellcolor[HTML]{34FF34}$+$       & \cellcolor[HTML]{FD6864}$-$       & \cellcolor[HTML]{34FF34}$+$       & \cellcolor[HTML]{FD6864}$-$       & \cellcolor[HTML]{FD6864}$-$       & \cellcolor[HTML]{34FF34}$+$       & \cellcolor[HTML]{FD6864}$-$       & \cellcolor[HTML]{FD6864}$-$       & \cellcolor[HTML]{FD6864}$-$       & \cellcolor[HTML]{FD6864}$-$       & \cellcolor[HTML]{FD6864}$-$       & \cellcolor[HTML]{FD6864}$-$       & \cellcolor[HTML]{FD6864}$-$       & \cellcolor[HTML]{FD6864}$-$       & \cellcolor[HTML]{FD6864}$-$       \\
$\LY$   & \cellcolor[HTML]{34FF34}$+$       & \cellcolor[HTML]{FD6864}$-$       & \cellcolor[HTML]{FD6864}$-$       & \cellcolor[HTML]{34FF34}$+$       & \cellcolor[HTML]{FD6864}$-$       & \cellcolor[HTML]{FD6864}$-$       & \cellcolor[HTML]{FD6864}$-$       & \cellcolor[HTML]{FD6864}$-$       & \cellcolor[HTML]{FD6864}$-$       & \cellcolor[HTML]{34FF34}$+$       & \cellcolor[HTML]{FD6864}$-$       & \cellcolor[HTML]{FD6864}$-$       & \cellcolor[HTML]{FD6864}$-$       & \cellcolor[HTML]{FD6864}$-$       & \cellcolor[HTML]{FD6864}$-$       & \cellcolor[HTML]{FD6864}$-$       \\
$\THO$  & \cellcolor[HTML]{FD6864}$-$       & \cellcolor[HTML]{FD6864}$-$       & \cellcolor[HTML]{FD6864}$-$       & \cellcolor[HTML]{FD6864}$-$       & \cellcolor[HTML]{FD6864}$-$       & \cellcolor[HTML]{FD6864}$-$       & \cellcolor[HTML]{FD6864}$-$       & \cellcolor[HTML]{FD6864}$-$       & \cellcolor[HTML]{FD6864}$-$       & \cellcolor[HTML]{FD6864}$-$       & \cellcolor[HTML]{FD6864}$-$       & \cellcolor[HTML]{FD6864}$-$       & \cellcolor[HTML]{FD6864}$-$       & \cellcolor[HTML]{FD6864}$-$       & \cellcolor[HTML]{FD6864}$-$       & \cellcolor[HTML]{FD6864}$-$       \\
$\FIWH$ & \cellcolor[HTML]{34FF34}$+$       & \cellcolor[HTML]{34FF34}$+$       & \cellcolor[HTML]{FD6864}$-$       & \cellcolor[HTML]{34FF34}$+$       & \cellcolor[HTML]{FD6864}$-$       & \cellcolor[HTML]{34FF34}$+$       & \cellcolor[HTML]{FD6864}$-$       & \cellcolor[HTML]{FD6864}$-$       & \cellcolor[HTML]{FD6864}$-$       & \cellcolor[HTML]{FD6864}$-$       & \cellcolor[HTML]{FD6864}$-$       & \cellcolor[HTML]{FD6864}$-$       & \cellcolor[HTML]{FD6864}$-$       & \cellcolor[HTML]{FD6864}$-$       & \cellcolor[HTML]{FD6864}$-$       & \cellcolor[HTML]{FD6864}$-$       \\
$\COO$  & \cellcolor[HTML]{34FF34}$+$       & \cellcolor[HTML]{34FF34}$+$       & \cellcolor[HTML]{FD6864}$-$       & \cellcolor[HTML]{34FF34}$+$       & \cellcolor[HTML]{34FF34}$+$       & \cellcolor[HTML]{34FF34}$+$       & \cellcolor[HTML]{34FF34}$+$       & \cellcolor[HTML]{FD6864}$-$       & \cellcolor[HTML]{34FF34}$+$       & \cellcolor[HTML]{34FF34}$+$       & \cellcolor[HTML]{FD6864}$-$       & \cellcolor[HTML]{FD6864}$-$       & \cellcolor[HTML]{34FF34}$+$       & \cellcolor[HTML]{FD6864}$-$       & \cellcolor[HTML]{34FF34}$+$       & \cellcolor[HTML]{34FF34}$+$       \\
$\JF$   & \cellcolor[HTML]{34FF34}$+$       & \cellcolor[HTML]{34FF34}$+$       & \cellcolor[HTML]{FD6864}$-$       & \cellcolor[HTML]{34FF34}$+$       & \cellcolor[HTML]{FD6864}$-$       & \cellcolor[HTML]{34FF34}$+$       & \cellcolor[HTML]{FD6864}$-$       & \cellcolor[HTML]{FD6864}$-$       & \cellcolor[HTML]{34FF34}$+$       & \cellcolor[HTML]{FD6864}$-$       & \cellcolor[HTML]{FD6864}$-$       & \cellcolor[HTML]{FD6864}$-$       & \cellcolor[HTML]{FD6864}$-$       & \cellcolor[HTML]{FD6864}$-$       & \cellcolor[HTML]{FD6864}$-$       & \cellcolor[HTML]{FD6864}$-$       \\
$\FIWF$ & \cellcolor[HTML]{34FF34}$+$       & \cellcolor[HTML]{34FF34}$+$       & \cellcolor[HTML]{FD6864}$-$       & \cellcolor[HTML]{34FF34}$+$       & \cellcolor[HTML]{FD6864}$-$       & \cellcolor[HTML]{34FF34}$+$       & \cellcolor[HTML]{FD6864}$-$       & \cellcolor[HTML]{FD6864}$-$       & \cellcolor[HTML]{34FF34}$+$       & \cellcolor[HTML]{FD6864}$-$       & \cellcolor[HTML]{34FF34}$+$       & \cellcolor[HTML]{FD6864}$-$       & \cellcolor[HTML]{FD6864}$-$       & \cellcolor[HTML]{FD6864}$-$       & \cellcolor[HTML]{FD6864}$-$       & \cellcolor[HTML]{FD6864}$-$       \\
$\B$    & \cellcolor[HTML]{34FF34}$+$       & \cellcolor[HTML]{34FF34}$+$       & \cellcolor[HTML]{FD6864}$-$       & \cellcolor[HTML]{34FF34}$+$       & \cellcolor[HTML]{FD6864}$-$       & \cellcolor[HTML]{34FF34}$+$       & \cellcolor[HTML]{34FF34}$+$       & \cellcolor[HTML]{FD6864}$-$       & \cellcolor[HTML]{FD6864}$-$       & \cellcolor[HTML]{34FF34}$+$       & \cellcolor[HTML]{FD6864}$-$       & \cellcolor[HTML]{FD6864}$-$       & \cellcolor[HTML]{FD6864}$-$       & \cellcolor[HTML]{FD6864}$-$       & \cellcolor[HTML]{FD6864}$-$       & \cellcolor[HTML]{34FF34}$+$       \\
$\MO$   & \cellcolor[HTML]{34FF34}$+$       & \cellcolor[HTML]{34FF34}$+$       & \cellcolor[HTML]{FD6864}$-$       & \cellcolor[HTML]{34FF34}$+$       & \cellcolor[HTML]{34FF34}$+$       & \cellcolor[HTML]{34FF34}$+$       & \cellcolor[HTML]{34FF34}$+$       & \cellcolor[HTML]{FD6864}$-$       & \cellcolor[HTML]{34FF34}$+$       & \cellcolor[HTML]{34FF34}$+$       & \cellcolor[HTML]{34FF34}$+$       & \cellcolor[HTML]{FD6864}$-$       & \cellcolor[HTML]{34FF34}$+$       & \cellcolor[HTML]{FD6864}$-$       & \cellcolor[HTML]{34FF34}$+$       & \cellcolor[HTML]{34FF34}$+$      
\end{tabular}
\end{center}
\end{figure}

\newpage

\begin{figure}[h!]
\begin{center}
\caption{Sporadic Groups Subquotient Table from the $\ATLAS$, II}
\label{table:Sporadic Groups Table from the ATLAS, II}
\begin{tabular}{ccccccccccc}
         & $\FIWW$                           & $\HN$                             & $\LY$                             & $\THO$                            & $\FIWH$                           & $\COO$                            & $\JF$                             & $\FIWF$                                         & $\B$                                            & $\MO$                                           \\
$\FIWW$ & \cellcolor[HTML]{FFFC9E}$\bullet$ & \cellcolor[HTML]{000000}          & \cellcolor[HTML]{000000}          & \cellcolor[HTML]{000000}          & \cellcolor[HTML]{000000}          & \cellcolor[HTML]{000000}          & \cellcolor[HTML]{000000}          & \cellcolor[HTML]{000000}                        & \cellcolor[HTML]{000000}                        & \cellcolor[HTML]{000000}                        \\
$\HN$   & \cellcolor[HTML]{FD6864}$-$       & \cellcolor[HTML]{FFFC9E}$\bullet$ & \cellcolor[HTML]{000000}          & \cellcolor[HTML]{000000}          & \cellcolor[HTML]{000000}          & \cellcolor[HTML]{000000}          & \cellcolor[HTML]{000000}          & \cellcolor[HTML]{000000}                        & \cellcolor[HTML]{000000}                        & \cellcolor[HTML]{000000}                        \\
$\LY$   & \cellcolor[HTML]{FD6864}$-$       & \cellcolor[HTML]{FD6864}$-$       & \cellcolor[HTML]{FFFC9E}$\bullet$ & \cellcolor[HTML]{000000}          & \cellcolor[HTML]{000000}          & \cellcolor[HTML]{000000}          & \cellcolor[HTML]{000000}          & \cellcolor[HTML]{000000}                        & \cellcolor[HTML]{000000}                        & \cellcolor[HTML]{000000}                        \\
$\THO$  & \cellcolor[HTML]{FD6864}$-$       & \cellcolor[HTML]{FD6864}$-$       & \cellcolor[HTML]{FD6864}$-$       & \cellcolor[HTML]{FFFC9E}$\bullet$ & \cellcolor[HTML]{000000}          & \cellcolor[HTML]{000000}          & \cellcolor[HTML]{000000}          & \cellcolor[HTML]{000000}                        & \cellcolor[HTML]{000000}                        & \cellcolor[HTML]{000000}                        \\
$\FIWH$ & \cellcolor[HTML]{34FF34}$+$       & \cellcolor[HTML]{FD6864}$-$       & \cellcolor[HTML]{FD6864}$-$       & \cellcolor[HTML]{FD6864}$-$       & \cellcolor[HTML]{FFFC9E}$\bullet$ & \cellcolor[HTML]{000000}          & \cellcolor[HTML]{000000}          & \cellcolor[HTML]{000000}                        & \cellcolor[HTML]{000000}                        & \cellcolor[HTML]{000000}                        \\
$\COO$  & \cellcolor[HTML]{FD6864}$-$       & \cellcolor[HTML]{FD6864}$-$       & \cellcolor[HTML]{FD6864}$-$       & \cellcolor[HTML]{FD6864}$-$       & \cellcolor[HTML]{FD6864}$-$       & \cellcolor[HTML]{FFFC9E}$\bullet$ & \cellcolor[HTML]{000000}          & \cellcolor[HTML]{000000}                        & \cellcolor[HTML]{000000}                        & \cellcolor[HTML]{000000}                        \\
$\JF$   & \cellcolor[HTML]{FD6864}$-$       & \cellcolor[HTML]{FD6864}$-$       & \cellcolor[HTML]{FD6864}$-$       & \cellcolor[HTML]{FD6864}$-$       & \cellcolor[HTML]{FD6864}$-$       & \cellcolor[HTML]{FD6864}$-$       & \cellcolor[HTML]{FFFC9E}$\bullet$ & \cellcolor[HTML]{000000}{\color[HTML]{000000} } & \cellcolor[HTML]{000000}{\color[HTML]{000000} } & \cellcolor[HTML]{000000}{\color[HTML]{000000} } \\
$\FIWF$ & \cellcolor[HTML]{34FF34}$+$       & \cellcolor[HTML]{FD6864}$-$       & \cellcolor[HTML]{FD6864}$-$       & \cellcolor[HTML]{FD6864}$-$       & \cellcolor[HTML]{34FF34}$+$       & \cellcolor[HTML]{FD6864}$-$       & \cellcolor[HTML]{FD6864}$-$       & \cellcolor[HTML]{FFFC9E}$\bullet$               & \cellcolor[HTML]{000000}                        & \cellcolor[HTML]{000000}                        \\
$\B$    & \cellcolor[HTML]{34FF34}$+$       & \cellcolor[HTML]{34FF34}$+$       & \cellcolor[HTML]{FD6864}$-$       & \cellcolor[HTML]{34FF34}$+$       & \cellcolor[HTML]{34FF34}$+$       & \cellcolor[HTML]{FD6864}$-$       & \cellcolor[HTML]{FD6864}$-$       & \cellcolor[HTML]{FD6864}$-$                     & \cellcolor[HTML]{FFFC9E}$\bullet$               & \cellcolor[HTML]{000000}                        \\
$\MO$   & \cellcolor[HTML]{34FF34}$+$       & \cellcolor[HTML]{34FF34}$+$       & \cellcolor[HTML]{FD6864}$-$       & \cellcolor[HTML]{34FF34}$+$       & \cellcolor[HTML]{34FF34}$+$       & \cellcolor[HTML]{34FF34}$+$       & \cellcolor[HTML]{FD6864}$-$       & \cellcolor[HTML]{34FF34}$+$                     & \cellcolor[HTML]{34FF34}$+$                     & \cellcolor[HTML]{FFFC9E}$\bullet$              
\end{tabular}
\end{center}
\end{figure}

\paragraph{Linear Representations, Projective Representations, and $X_G$}\label{par:Linear Representations, Projective Representations, and XG}
Let $G$ be a finite simple group. Every non-trivial linear and projective representation of $G$ is faithful, so the corresponding quotient maps yield upper bounds on $\RD_k(G)$ \cite[Example 4.6]{FarbKisinWolfson2023}. For each sporadic group $G$, we now compare $X_G$ with a minimal-dimensional linear representation $W_G$ and a minimal-dimensional projective representation $\PP\lp V_G \rp$ (continuing with the notation of Notation \ref{not:Representations of Sporadic Groups}). In the cases where the linear representation comes from the Schur cover, we include the Schur cover as well.

\begin{figure}[h!]
\begin{center}
\caption{Minimal Linear Representations, Projective Representations, and $X_G$}
\label{table:Dimensions of Linear Representations, Projective Representations, and XG}
\begin{tabular}{|c|c|c|c|}
\hline
Group $G$ &$\dim\lp W_G \rp$ &$\dim\lp \PP\lp V_G \rp \rp$ &$\dim\lp X_G \rp$\\
\hline
$\MOO$ &10 &9 &6 \\
\hline
$\MOW$ &11 &10 &7 \\
\hline
$\MWH$ &22 &21 &17 \\
\hline
$\HS$ &22 &21 &18 \\
\hline
$\MCL$ &22 &21 &19 \\
\hline
$\MWF$ &23 &22 &18 \\
\hline
$\COH$ &23 &22 &20 \\
\hline
$\COW$ &23 &22 &20 \\
\hline
$\HE$ &51 &50 &48 \\
\hline
$\JO$ &56 &55 &51 \\
\hline
$\FIWW$ &78 &77 &74 \\
\hline
$\HN$ &133 &132 &129 \\
\hline
$\THO$ &248 &247 &244 \\
\hline
$\FIWH$ &782 &781 &776\\
\hline
$\JF$ &1333 &1332 &1328 \\
\hline
$\LY$ &2480 &2479 &2475 \\
\hline
$\B$ &4371 &4370 &4365 \\
\hline
$\MO$ &196883 &196882 &196874 \\
\hline
\end{tabular}
\end{center}
\end{figure}

\begin{figure}[h!]
\begin{center}
\caption{Minimal Linear Representations, Projective Representations, $X_G$, and the Schur Covers}
\label{table:Dimensions of Linear Representations, Projective Representations, XG, and the Schur Multipliers}
\begin{tabular}{|c|c|c|c|c|}
\hline
Group $G$ &$\dim\lp W_G \rp$ &$\dim\lp \PP\lp V_G \rp \rp$ &$\dim\lp X_G \rp$ &$a.G$\\
\hline
$\JW$ &14 &5 &5 &$2.\JW$ \\
\hline
$\MWW$ &20 &9 &8 &$12.\MWW$ \\
\hline
$\SUZ$ &143 &11 &10 &$6.\SUZ$ \\
\hline
$\JH$ &85 &17 &16 &$3.\JH$ \\
\hline
$\COO$ &276 &23 &21 &$2.\COO$ \\
\hline
$\RU$ &378 &27 &26 &$2.\RU$ \\
\hline
$\ON$ &10944 &341 &338 &$3.\ON$ \\
\hline
$\FIWF$ &8671 &782 &779 &$3.\FIWF$\\
\hline
\end{tabular}
\end{center}
\end{figure}

\newpage

\appendix

\section{Power Series Expansions of Molien Series}\label{app:Power Series Expanions of Molien Series}

For each sporadic group $G$, we record the initial terms of the power series expansions of $M\lp \rho_G;t \rp$ in Figure \ref{table:Power Series Expansions}.

\begin{figure}[h!]
\begin{center}
\caption{Power Series Expansions}
\label{table:Power Series Expansions}
\begin{tabular}{|c|c|}
\hline
Group &Initial Terms for the Power Series Expansion of $M\lp \rho_G; t \rp$\\
\hline
$\JW$ &$1 + t^{12} + t^{20} + 2t^{24} + t^{28} + 2t^{30} + 3t^{32} + t^{34} + 4t^{36} + 2t^{38} + O\lp t^{40} \rp$ \\
\hline
$\MOO$ &$1 + t^{2} + t^{3} + 2t^{4} + 3t^{5} + 5t^{6} + 6t^{7} + 11t^{8} + 16t^{9} + 26t^{10} + 38t^{11} + 61t^{12} + 91t^{13} + O\lp t^{14} \rp$\\
\hline
$\MOW$ &$1 + t^{2} + t^{3} + 2t^{4} + 2t^{5} + 5t^{6} + 4t^{7} + 9t^{8} + 10t^{9} + 17t^{10} + 20t^{11} + 36t^{12} + 39t^{13} + 67t^{14} + O\lp t^{15} \rp$ \\
\hline
$\MWW$ & $1 + t^4 + t^6 + 2t^8 + 3t^{10} + 6t^{12} + 9t^{14} + 15t^{16} + 26t^{18} + O\lp t^{20} \rp$\\
\hline
$\SUZ$ &$1 + t^{12} + t^{18} + 3t^{24} + 3t^{30} + 7t^{36} + O\lp t^{40} \rp$ \\
\hline
$\JH$ &$1 + t^{6} + t^{9} + 10t^{12} + 26t^{15} + 143t^{18} + 680t^{21} + 3310t^{24} + 14229t^{27} + 55826t^{30} + O\lp t^{33} \rp$\\
\hline
$\MWH$ &$1 + t^{2} + t^{3} + 2t^{4} + 3t^{5} + 6t^{6} + 9t^{7} + 17t^{8} + 27t^{9} + 49t^{10} + 86t^{11} + 159t^{12} + 292t^{13} + O\lp t^{14} \rp$ \\
\hline
$\HS$ &$1 + t^{2} + 2t^{4} + t^{5} + 5t^{6} + 3t^{7} + 12t^{8} + 9t^{9} + 29t^{10} + 28t^{11} + 77t^{12} + 87t^{13} + 220t^{14} + O\lp t^{15} \rp $\\
\hline
$\MCL$ &$1 + t^{2} + t^{4} + t^{5} + 2t^{6} + 3t^{7} + 5t^{8} + 6t^{9} + 10t^{10} + 14t^{11} + 21t^{12} + 29t^{13} + 48t^{14} + 70t^{15} + O\lp t^{16} \rp $\\
\hline
\end{tabular}
\end{center}
\end{figure}

\begin{center}
\begin{tabular}{|c|c|}
\hline
$\MWF$ &$1 + t^{2} + t^{3} + 2t^{4} + 2t^{5} + 5t^{6} + 5t^{7} + 11t^{8} + 14t^{9} + 25t^{10} + 35t^{11} + 65t^{12} + 89t^{13} + O\lp t^{14} \rp$\\
\hline
$\COH$ &$1 + t^{2} + t^{4} + 2t^{6} + 3t^{8} + t^{9} + 5t^{10} + 2t^{11} + 9t^{12} + 3t^{13} + 14t^{14} + 7t^{15} + 23t^{16} + 13t^{17} + O\lp t^{18} \rp $\\
\hline
$\COW$ &$1 + t^{2} + t^{4} + t^{6} + 2t^{8} + 3t^{10} + t^{11} + 5t^{12} + t^{13} + 7t^{14} + 2t^{15} + 11t^{16} + 3t^{17} + 16t^{18} + O\lp t^{19} \rp$\\
\hline
$\COO$ &$1 + t^{2} + t^{4} + t^{6} + t^{8} + t^{10} + 2t^{12} + 2t^{14} + 3t^{16} + 4t^{18} + O\lp t^{20} \rp$ \\
\hline
$\RU$ &$1 + t^{4} + 2t^{8} + 6t^{12} + 2t^{14} + 22t^{16} + 27t^{18} + 154t^{20} + 439t^{22} + 1966t^{24} + 7189t^{26} + O\lp t^{28} \rp $\\
\hline
$\HE$ &$1 + t^{3} + t^{4} + t^{5} + 4t^{6} + 5t^{7} + 13t^{8} + 30t^{9} + 82t^{10} + 245t^{11} + 907t^{12} + 3424t^{13} + O\lp t^{14} \rp $\\
\hline
$\JO$ &$1 + t^{2} + t^{3} + 8t^{4} + 34t^{5} + 361t^{6} + 2820t^{7} + 22346t^{8} + 156939t^{9} + 1021469t^{10} + O\lp t^{11} \rp $\\
\hline
$\FIWW$ &$1 + t^2 + t^4 + 2t^6 + 5t^8 + t^9 + 13t^{10} + 4t^{11} + 60t^{12} + 31t^{13} + 488t^{14} + 912t^{15} + O\lp t^{16} \rp $ \\
\hline
$\HN$ &$1 + t^2 + t^4 + 2t^6 + t^7 + 5t^8 + 6t^9 + 27t^{10} + 92t^{11} + 637t^{12} + 5018t^{13} + 47239t^{14} + O\lp t^{15} \rp $\\
\hline
$\THO$ &$1 + t^2 + t^4 + t^6 + 4t^8 + 15t^{10} + 50t^{11} + 1854t^{12} + 31610t^{13} + 607473t^{14} + O\lp t^{15} \rp$  \\
\hline
$\ON$ &$1 + 16t^6 + 426595t^9 + 14039408007t^{12} + 230067642077481t^{15} + O\lp t^{18} \rp $ \\
\hline
$\FIWH$ &$1 + t^2 + t^3 + 2t^4 + 3t^5 + 9t^6 + 15t^7 + 57t^{8} + 324t^9 + 7961t^{10} + 456255t^{11} + O\lp t^{12} \rp$ \\
\hline
$\FIWF$ &$1 + t^3 + 3t^6 + 11t^9 + 355t^{12} + 17843536t^{15} + 1848868683076t^{18} + O\lp t^{21} \rp$ \\
\hline
$\JF$ &$1 + t^4 + 2t^6 + 2t^7 + 31t^8 + 521t^9 + 60960t^{10} + 7118797t^{11} + 795955946t^{12} + O\lp t^{13} \rp$ \\
\hline
$\LY$ &$1 + 23t^6 + 21041t^7 + 697156t^8 + 191631120t^9 + 47708455027t^{10} + O\lp t^{11} \rp$ \\
\hline
$\B$ &$1 + t^2 + 2t^4 + 3t^6 + 7t^8 + 20t^{10} + 3t^{11} + 243t^{12} + 8164t^{13} + 2665262t^{14} + O\lp t^{15} \rp$ \\
\hline
$\MO$ &$1 + t^2 + t^3 + 2t^4 + 2t^5 + 6t^6 + 6t^7 + 16t^8 + 27t^9 + 68t^{10} + 182t^{11} + 956t^{12} + O\lp t^{13} \rp $\\
\hline
\end{tabular}
\end{center}



\paragraph{Author:} Claudio G\'{o}mez-Gonz\'{a}les
\vspace{-6pt}
\paragraph{Affiliation:} Carleton College; Department of Mathematics and Statistics
\vspace{-6pt}
\paragraph{E-mail:} cgonzales@carleton.edu

\vspace{16pt}

\paragraph{Author:} Alexander J. Sutherland
\vspace{-6pt}
\paragraph{Affiliation:} The Ohio State University; Department of Mathematics
\vspace{-6pt}
\paragraph{E-mail:} sutherland.159@osu.edu

\vspace{16pt}

\paragraph{Author:} Jesse Wolfson
\vspace{-6pt}
\paragraph{Affiliation:} University of California, Irvine; Department of Mathematics
\vspace{-6pt}
\paragraph{E-mail:} wolfson@uci.edu

\vspace{16pt}

\paragraph{MSC Classes:} 14L30 (Primary); 13A50, 20C25, 20C34 (Secondary)

\paragraph{Keywords:} resolvent degree, torsors, versality, rational points, sporadic groups, Monster group

\end{document}